\documentclass[3p]{elsarticle}

\usepackage{hyperref}

\journal{Applied Mathematics and Computation}

\usepackage{amsopn,amsfonts}
\usepackage{multirow}
\usepackage{algorithm,algorithmic}

\DeclareMathOperator{\R}{\mathbb{R}}
\DeclareMathOperator{\N}{\mathbb{N}}
\DeclareMathOperator{\M}{\mathcal{M}}
\DeclareMathOperator{\argmin}{argmin}
\newcommand{\IntEnt}[1]{\left[\!\left[#1\right]\!\right]}

\newtheorem{thm}{Theorem}
\newtheorem{prop}[thm]{Proposition}
\newdefinition{df}{Definition}
\newdefinition{rmk}{Remark}
\newproof{pf}{Proof}

\usepackage[breakable,theorems]{tcolorbox}
\newtcbtheorem[auto counter]{Dcp}{DC Programming Formulation}{coltitle=black,colbacktitle=gray!20,colframe=black,colback=white,coltext=black}{Dcp}

\begin{document}

\begin{frontmatter}

\title{Improved DC Programming Approaches for Solving the Quadratic Eigenvalue Complementarity Problem \footnote{The paper is accepted by Applied Mathematics and Computation. \url{https://doi.org/10.1016/j.amc.2019.02.017}}}

\author[mainaddress,mainaddressbis]{Yi-Shuai Niu\corref{correspondingauthor}}
\address[mainaddress]{School of Mathematical Sciences, Shanghai Jiao Tong University, Shanghai 200240, China}
\address[mainaddressbis]{SJTU-Paristech Elite Institute of Technology, Shanghai Jiao Tong University, Shanghai 200240, China}
\cortext[correspondingauthor]{Corresponding author}
\ead{niuyishuai@sjtu.edu.cn}

\author[secondaryaddress]{Joaquim J\'udice}
\address[secondaryaddress]{Instituto de Telecomunica\c{c}\~{o}es, Coimbra, Portugal}

\author[thirdaddress]{Hoai An Le Thi}
\address[thirdaddress]{University of Lorraine, Metz, France}

\author[fourthaddress]{Dinh Tao Pham}
\address[fourthaddress]{National Institute of Applied Sciences, Rouen, France}

\begin{abstract}
In this paper, we discuss the solution of a Quadratic Eigenvalue Complementarity Problem (QEiCP) by using Difference of Convex (DC) programming approaches. We first show that QEiCP can be represented as dc programming problem. Then we investigate different dc programming formulations of QEiCP and discuss their dc algorithms based on a well-known method -- DCA. A new local dc decomposition is proposed which aims at constructing a better dc decomposition regarding to the specific feature of the target problem in some neighborhoods of the iterates. This new procedure yields faster convergence and better precision of the computed solution. Numerical results illustrate the efficiency of the new dc algorithms in practice. 
\end{abstract}


\end{frontmatter}


\section{Introduction}
Let $\M_n$ be the set of all $n\times n$ real matrices with $n\in \N^*$, where $\N^*$ denotes the set of all nonzero natural numbers. Given three matrices $A,B,C$ in $\M_n$, the \emph{Quadratic Eigenvalue Complementarity Problem} (QEiCP) consists of finding a scalar $\lambda \in \mathbb{R}$ and a vector $x\in \mathbb{R}^n\setminus\{0\}$ such that
\begin{equation}\label{prob:QEiCP}
\left\lbrace \begin{array}{l}
w=\lambda^2Ax+\lambda Bx+Cx,\\
x^{\top}w=0,\\
x\geq 0, w\geq 0.
\end{array}\right. 
\end{equation}
This problem has been introduced in \cite{Seeger11} where some applications were highlighted. We denote the problem (\ref{prob:QEiCP}) as QEiCP($A,B,C$). In any solution $(\lambda,x)$ of QEiCP, the $\lambda$-component is a \emph{quadratic complementary eigenvalue} and the $x$-component is a \emph{quadratic complementary eigenvector} of $(A,B,C)$ associated to $\lambda$.

QEiCP is an extension of the well-known Eigenvalue Complementarity Problem (EiCP). This last problem has been introduced in \cite{Seeger99} and consists of finding a complementary eigenvalue $\lambda\in \R$ and an associated complementary eigenvector $x\in \R^n\setminus\{0\}$ such that 
\begin{equation}\label{prob:EiCP}
\left\lbrace \begin{array}{l}
w=\lambda Bx-Cx,\\
x\geq 0, w\geq 0,\\
x^{\top}w=0.\\
\end{array}\right.
\end{equation}  
where $B,C\in \M_n$. Therefore, EiCP is a special case of QEiCP when the matrix $A$ is null (i.e., EiCP($B,C$) = QEiCP($0,B,-C$)). 

EiCP has a solution if the matrix $B$ of the leading $\lambda$-term is positive definite (PD) \cite{Seeger99,Judice09}, that is,
\begin{equation}
x^{\top} B x >0, \forall x\in \R^n\setminus\{0\}.
\end{equation} 
During the past several years, many theoretical results \cite{Seeger99,Judice09,Bras12,Costa10}, applications \cite{Costa04,Fernandes17,Patrascu15} and extensions \cite{Adly15,Chen16,Costa17,Fernandes16,Zhang18,Yihui09} of EiCP have been discussed and a number of efficient algorithms have been proposed for the solution of this problem \cite{Adly11,Bras17,Facchinei,Fernandes14,Iusem_inpress,Judice08,Lethi12,Niu12}. 

Contrary to the EiCP, QEiCP may have no solution even when the matrix $A$ of the leading $\lambda$-term is PD \cite{Bras16}. The existence of a solution of QEiCP depends on the properties of the given matrices $A,B$ and $C$. If the matrix $A$ is PD, then QEiCP has a solution if one of the two following conditions holds:
\begin{itemize}[(i)]
	\item $C\notin S_0$ \cite{Bras16}, where $S_0$ is the class of matrices defined in \cite{Cottle} by
	\begin{equation}
	S_0=\{C\in \M_n: \exists x\in \R^n_+\setminus\{0\}, Cx\geq 0\}.
	\end{equation}
	\item \emph{co-hyperbolic} \cite{Seeger11} \begin{equation}\label{eq:co-hyperbolic}
	(x^{\top}Bx)^2\geq 4(x^{\top}Ax)(x^{\top}Cx), \forall x\geq 0.
	\end{equation}
\end{itemize}
As discussed in \cite{Bras16}, none of these conditions implies the other. Furthermore,  investigating whether $C\notin S_0$ reduces to solving a special linear program \cite{Bras16}. On the other hand, it is very hard to prove that the co-hyperbolic condition holds in practice. However, there are some sufficient conditions which imply (\ref{eq:co-hyperbolic}). For instance, this occurs if $A$ and $-C$ are both PD matrices.

A number of algorithms have been proposed for the solution of QEiCP and its extensions to other cases when $A\in$ PD and one of the conditions $C\notin S_0$ or co-hyperbolicity holds \cite{Adly11,Bras15,Bras16,Fernandes,Fernandes14,Iusem17,Iusem16,Seeger11}. Two nonlinear programming (NLP) formulations of QEiCP have been introduced in \cite{Bras16,Fernandes14} such that $(\lambda,x)$ is a solution of QEiCP if and only if $(\lambda,x)$ is a global minimum of NLP with an optimal value equal to zero. In this paper, we introduce dc programming formulations for QEiCP which are based on these nonlinear programs.

As for the NLPs mentioned before, a stationary point of the dc program with null objective function value gives a solution of QEiCP. In order to find such a stationary point, the classical DC algorithm (DCA) \cite{Pham98,Pham18} may be used. In our previous work \cite{Niu15}, we showed that DCA has a very good performance for solving QEiCP with fast convergence to a complementary eigenvalue for most of the test problems. However, it is also observed that for some instances, DCA converges quite slowly with a very small reduction of the objective function value from one iteration to the next. Improving the performance of DCA for QEiCP is the biggest motivation for this paper. Next, we summarize the major contributions of this paper.
\begin{enumerate}[(i)]
	\item Finding an appropriate dc decomposition for any dc function is an open question in dc programming \cite{Pham18}. We propose a tool to \emph{qualify an appropriate dc decomposition} (Proposition \ref{prop:betterdcd}) and Definition \ref{def:betterdcp}) which relies on the convexity of the dc components. Since any dc function has infinitely many dc decompositions, this tool is important for comparing the quality of two dc decompositions.   
	\item We propose a new dc programming formulation for QEiCP based on the \emph{Difference of Convex Sums-of-Squares (DCSOS) decompositions} defined in \cite{Niu17}. This new decomposition is compared with the classical universal dc decomposition proposed in \cite{Niu15}. Numerical results reported in section \ref{sec:experiments} and illustrated in Tables \ref{tab:numresultsdcas} and \ref{tab:numresultsdcasbis} indicate a good performance of our new approach.  
	\item We introduce a new formulation to \emph{compute a tighter bound} for complementary eigenvalues given in Theorem \ref{thm:boundsoflambdabis}. Numerical simulations illustrate that this bound is tighter than those given in \cite{Fernandes14,Niu15}. A tighter bound is still required for the universal dc decomposition, which will provide a better dc decomposition for DCA.
	\item The use of DCA for solving DCSOS decomposition requires solving a convex polynomial optimization in each iteration that is observed inefficient particularly for large-scale instances using existing convex optimization solvers such as IPOPT \cite{Ipopt}. In order to improve the efficiency of DCA, we propose an \emph{equivalent second order conic programming (SOCP) formulation} for this convex polynomial optimization problem. The SOCP formulation can be solved more effectively in practice by existing SOCP solvers such as CPLEX \cite{Cplex} and GUROBI \cite{Gurobi}.  
	\item Finding a good initial point is an open question in dc programming \cite{Pham18}. We propose a heuristic procedure to \emph{estimate a good initial point} for DCA when applied to QEiCP. This estimation only requires solving a strictly convex quadratic optimization problem over a standard simplex.
	\item We propose a \emph{local dc decomposition algorithm} which updates in each iteration a better dc decomposition and yields a faster convergence rate for DCA to get a computed solution with a better precision.
\end{enumerate}

The organization of the paper is as follows. In section \ref{sec:betterdcdecomp}, we discuss how to qualify the dc decompositions which provide an useful tool for analyzing the quality of different dc decompositions. Section \ref{sec:NLPformulations} contains two different nonlinear programming formulations (NLPs) for QEiCP and their dc programming formulations are established in sections \ref{sec:dcpofP}, \ref{sec:dcpofPhat} and \ref{sec:dcpofP'}. Some new results concerning tighter bounds estimation for complementary eigenvalues are discussed in section \ref{sec:boundoflambda}. Some algorithms based on the classical DCA for solving the dc programming formulations are presented in section \ref{sec:dcas}. In particular, the SOCP formulation for polynomial convex subproblem is introduced in subsection \ref{subsec:socpformulation}, a heuristic for a good initial point estimation is discussed in subsection \ref{subsec:startingpoint}, and a new local dc decomposition method is proposed in subsection \ref{subsec:localdcdecomp}. Computational experiments with these algorithms are reported in section \ref{sec:experiments}. Finally, some conclusions are presented in the last section of the paper.

\section{Quality of dc decompositions}\label{sec:betterdcdecomp}
Let $D$ be a non-empty closed convex set of $\R^n$, let $f:D \to \R$ be a dc function with a dc decomposition as $f=g-h$ where $g$ and $h$ are real valued convex functions defined on $D$. A dc program is defined by \begin{equation}
\min\{g(x)-h(x):x\in D\}
\end{equation}
An efficient DC algorithm, called DCA, was introduced by D.T. Pham in 1985 and extensively developed by H.A. Le Thi and D.T. Pham since 1994 (the reader can refer to  \url{http://www.lita.univ-lorraine.fr/~lethi/index.php} and \cite{Pham98,Pham05,Pham18}). 

DCA consists of constructing two sequences $\{x^k\}$ and $\{y^k\}$ via the scheme:
\begin{equation*}
\begin{array}{c}
\text{Given } x^0\in \R^n,\\
x^{k}\rightarrow y^{k}\in \partial
h(x^{k})\\
\swarrow\\
x^{k+1}\in \partial g^{*}(y^{k})=\argmin\{g(x)-\langle x,y^k\rangle: x\in D\}.
\end{array}
\end{equation*}

The symbol $\partial h$ stands for the sub-differential of the convex function $h$, and $g^*$ is the conjugate function of $g$. These definitions are fundamental to understand the algorithm and can be found in any textbook of convex analysis (see e.g., \cite{Rockafellar}).

There exist infinitely many dc decompositions for any dc function $f$. In fact, if $f$ has a dc decomposition $g-h$, then $f = (g+\phi) - (h+\phi)$ is also a dc decomposition for any convex function $\phi:D\to \R$. An open question in the field of dc programming is what is a good dc decomposition for DCA and how to find it. 

Firstly, we should define some rules to qualify a dc decomposition. Regarding to the geometrical interpretation of DCA given in \cite{Niu10,Niu14}, the quality of a dc decomposition depends on the specific feature of the dc function $f$ at a given point in the convex set $D$. More specifically, 
DCA requires constructing in the $k$-th ($k\in \N^*$) iteration a convex overestimation of $f$ at point $x^k$, denoted by $f^k$ as \begin{equation}
f^k(x) = g(x)-(h(x^{k})+\langle x-x^{k},y^{k}\rangle)
\end{equation} 
where $y^k\in \partial h(x^k)$. So $f^k(x)$ is obtained by replacing $h(x)$ by its affine approximation at $x^k$. Hence  
\begin{equation}\label{eq:fvsfk}
f^k(x)\geq f(x), \forall x\in D, \forall k\in \N^*.
\end{equation} 

Therefore, the smaller the gap between $f^k$ and $f$ on $D$ is, the better $f^k$ fits $f$, and the better minimizers of $f^k$ approach minimizers of $f$. 

In order to find out what is the most important key point to reduce the gap between the convex function $f^k$ and the dc function $f$, we recall that for two given convex functions $g$ and $g'$, we say that $g$ is more (resp. less) convex than $g'$ if $g-g'$ (resp. $g'-g$) is still a convex function. Then, we have the following result:
\begin{prop}\label{prop:betterdcd}
	Let $g-h$ and $g'-h'$ be two dc decompositions for a dc function $f$ over a closed convex set $D$. 
	\begin{itemize}
		\item[(i)] The function $g$ is less (resp. more) convex than $g'$ if and only if $h$ is less (resp. more) convex than $h'$.
		\item[(ii)] If $g$ is less convex than $g'$, then the convex overestimations $f_g^k$ and $f_{g'}^k$ given by  \[f_g^k(x) = g(x) - (h(x^{k})+\langle x-x^{k},y^{k}\rangle), \text{ with } y^k\in\partial h(x^k). \]
		\[f_{g'}^k(x) = g'(x) - (h'(x^{k})+\langle x-x^{k},y'^{k}\rangle), \text{ with } y'^k\in\partial h'(x^k).\]
		satisfy \[f_g^k(x)\leq f_{g'}^k(x),\forall x\in D.\]
	\end{itemize}
\end{prop}
\begin{pf}
	(i) If $g$ is less (resp. more) convex than $g'$, then $g'-g$ (resp. $g-g'$) is a convex function, and $h'-h = g'- g$ (resp. $h-h' = g-g'$) is also convex, which yields the desired result.
	
	\noindent(ii) Let us denote $d = g'-g$, which is a convex function on $D$ if $g$ is less convex than $g'$. Then $g = g'-d$ and $h=h'-d$. Thus we have
	\begin{equation}\nonumber
	\begin{array}{rl}
	f_g^k(x) &= g(x) - (h(x^{k})+\langle x-x^{k},y^{k}\rangle) \\
	&= g'(x)-d(x) - (h'(x^{k}) - d(x^k) + \langle x-x^{k},y^{k}-y'^k+y'^k\rangle)\\
	&= g'(x) - (h'(x^k) + \langle x-x^k,y'^k\rangle) - (d(x) - (d(x^k) + \langle x-x^k,y'^k-y^k\rangle))\\
	&= f_{g'}^k(x) - (d(x) - (d(x^k) + \langle x-x^k,y'^k-y^k\rangle))
	\end{array}
	\end{equation}
	
	Since $d$ is convex and $y'^k-y^k \in \partial (h'-h) (x^k) = \partial d(x^k)$, then we get \[d(x) - (d(x^k) + \langle x-x^k,y'^k-y^k\rangle) \geq 0, \forall x\in D.\]
	
	Hence, \[f_g^k(x)\leq f_{g'}^k(x),\forall x\in D.\]
\end{pf}

According to Proposition \ref{prop:betterdcd} and the inequality (\ref{eq:fvsfk}), if $g$ is less convex than $g'$ then \[f(x)\leq f_g^k(x)\leq f_{g'}^k(x),\forall x\in D,\]
which means that $f_g^k$ should fit better $f$ than $f_{g'}^k$ on $D$. This result reveals that the most important point for choosing a good dc decomposition is to use functions $g$ and $h$ that are less convex as possible. Now, we can give a definition that states when dc decomposition is better than another one.
\begin{df}\label{def:betterdcp}
	Let $g-h$ and $g'-h'$ be two dc decomposition for a dc function $f$. We say that the dc decomposition with $g$ and $h$ is better than the one with $g'$ and $h'$ if $g$ is less convex than $g'$ or $h$ is less convex than $h'$. 
\end{df}

\section{Nonlinear programming formulations of QEiCP}\label{sec:NLPformulations}
The QEiCP defined in (\ref{prob:QEiCP}) is obviously equivalent to the problem:
\begin{equation}\label{prob:QEiCPbis}
\left\lbrace \begin{array}{l}
w=\lambda^2Ax+\lambda Bx+Cx,\\
x^{\top}w=0,\\
e^{\top}x=1,\\
x\geq 0, w\geq 0.
\end{array}\right. 
\end{equation}
by introducing the additional constraint $e^{\top}x = 1$. As previously discussed in  \cite{Fernandes14,Niu15}, the problem (\ref{prob:QEiCPbis}) is equivalent to :
\begin{equation}\label{prob:QEiCPbisbis}
\left\lbrace \begin{array}{l}
w=Az+By+Cx,\\
y=\lambda x,\\
z=\lambda y,\\
x^{\top}w=0,\\
e^{\top}x=1,\\
e^{\top}y=\lambda,\\
x\geq 0, w\geq 0.
\end{array}\right. 
\end{equation}
which can be rewritten as a nonlinear polynomial optimization problem described as follows: 
\begin{equation}\label{prob:nlp}
(P) \qquad 
\min \{f(x,y,z,w,\lambda): (x,y,z,w,\lambda) \in \mathcal{C} \}\nonumber
\end{equation}
where \begin{equation}\label{eq:f}
f(x,y,z,w,\lambda) = \|y-\lambda x\|^2 + \|z-\lambda y\|^2 + x^{\top}w
\end{equation} and \begin{equation}\label{eq:C}
\mathcal{C} = \{(x,y,z,w,\lambda): w=Az+By+Cx, e^{\top}x=1, e^{\top}y=\lambda, (x,w,z)\geq 0\}. 
\end{equation}
In the whole paper, we denote $\|x\|$ for the euclidean norm (2-norm) of vector $x\in \R^n$.

\begin{rmk}
	In problem (\ref{prob:QEiCPbisbis}), the linear constraint $e^{\top}y=\lambda$ is redundant since it is always satisfied when $y=\lambda x$ and $e^{\top}x=1$. However, it is useful to be presented in $(P)$ when we remove $y=\lambda x$ and add $\|y-\lambda x\|$ into objective function, since this linear constraint reveals a relationship between $y$ and $\lambda$, which could reduce the feasible set of $(P)$ without loosing any optimal solution. 
\end{rmk}

The following theorem establishes the equivalence between QEiCP and $(P)$ which is a simple consequence of the definition of QEiCP and $(P)$.
\begin{thm}\label{thm:equivofQEiCPandP} \mbox{}\par
	\begin{enumerate}
		\item[(i)] If $(P)$ is infeasible, then QEiCP($A,B,C$) is infeasible.
		\item[(ii)] If $(P)$ has nonzero global optimal value, then QEiCP($A,B,C$) has no solution.  
		\item[(iii)] If $(P)$ has a global optimal solution $(x^*,y^*,z^*,w^*,\lambda^*)$ such that $f(x^*,y^*,z^*,w^*,\lambda^*)=0$, then $(\lambda^*,x^*)$ is a solution of QEiCP($A,B,C$). 
		\item[(iv)] Conversely, if $(\lambda^*,x^*)$ is a solution of QEiCP($A,B,C$), then $(\bar{x},\bar{y},\bar{z},\bar{w},\bar{
			\lambda})$ defined by
		$$\bar{\lambda}=\lambda^*; \bar{x} = \frac{x^*}{e^{\top}x^*}; \bar{y}= \bar{\lambda}\bar{x}; \bar{z}=\bar{\lambda}\bar{y}; \bar{w} = A\bar{z} + B\bar{y} + C\bar{x}$$
		is a global optimal solution of $(P)$ with zero global optimal value.
	\end{enumerate}
\end{thm}

Observe that under the assumption $x\geq 0, w\geq 0$, the complementarity constraint $w^{\top}x=0$ is equivalent to $\sum_{i=1}^{n} \min\{x_i,w_i\} = 0$. Hence the bilinear term $w^{\top}x$ in $f$ can be replaced by a concave polyhedral function $\sum_{i=1}^{n} \min\{x_i,w_i\}$ which yields the following equivalent formulation of $(P)$ 
\begin{equation}\label{prob:nlp'}
(P') \qquad 
\min \{ f'(x,y,z,w,\lambda) : (x,y,z,w,\lambda)\in \mathcal{C} \}\nonumber
\end{equation}
where \begin{equation}
f'(x,y,z,w,\lambda)= \|y-\lambda x\|^2 + \|z-\lambda y\|^2 + \sum_{i=1}^{n} \min\{x_i,w_i\}
\end{equation} and $\mathcal{C}$ is defined in (\ref{eq:C}).
The difference between $(P)$ and $(P')$ is only related to their objective functions. The objective function of the problem $(P)$ is a $\mathcal{C}^{\infty}$ polynomial function, while the objective function in $(P')$ consists of a polyhedral function which is not differentiable at some points. Nevertheless the non-smoothness of problem $(P')$, it is worthwhile to introduce this formulation since this polyhedral function can be easily presented as a DC function and DCA possesses finite convergence for optimizing a polyhedral DC function on a convex set defined by linear constraints \cite{Pham02}. In the following sections, we will investigate how to represent the problems $(P)$ and $(P')$ into dc programming formulations.

\section{A first dc programming formulation}\label{sec:dcpofP}
Problem $(P)$ is a nonconvex polynomial optimization problem in which a nonconvex polynomial function $f$ is minimized on a polyhedral convex set $\mathcal{C}$. As in \cite{Niu15}, we reformulate $(P)$ into a dc programming problem by writing the function $f$ given by (\ref{eq:f}) as follows:
\begin{align}\label{eq:df}
f(x,y,z,w,\lambda) & =\|y\|^2 + \|z\|^2 -2\lambda y^{\top}(x+z) + \lambda^2 (\|x\|^2 + \|y\|^2) + x^{\top}w \\
& =f_0(y,z) + f_1(x,y,z,\lambda) + f_2(x,y,\lambda) + f_3(x,w) \nonumber
\end{align}
with 
\begin{equation}\label{eq:f0-f3}
\left\{\begin{array}{l}
f_0(y,z)=\|y\|^2 + \|z\|^2,\\
f_1(x,y,z,\lambda)=-2\lambda y^{\top}(x+z),\\
f_2(x,y,\lambda)= \lambda^2 (\|x\|^2 + \|y\|^2),\\ 
f_3(x,w)= x^{\top}w.
\end{array}
\right.
\end{equation}

The function $f_0$ is a convex quadratic function, while $f_1,f_2$ and $f_3$ are all nonconvex polynomial functions. Next, we explain how to get a dc decomposition for $f_1,f_2$ and $f_3$.

\begin{enumerate}
	\item As in \cite{Niu15}, a dc decomposition for bilinear function $f_3$ is given by
	\begin{equation}\label{eq:dcdf3}
	f_3(x,w)= \frac{\|x+w\|^2}{4}-\frac{\|x-w\|^2}{4}
	\end{equation}
	where $\frac{\|x+w\|^2}{4}$ and $\frac{\|x-w\|^2}{4}$ are both convex quadratic functions.
	
	\item A dc decomposition for $f_2$ is given by
	\begin{equation}\label{eq:dcdf2}
	f_2(x,y,\lambda) = \frac{(\lambda^2 + \|x\|^2)^2 + (\lambda^2 + \|y\|^2)^2}{2} - \frac{2\lambda^4 + \|x\|^4 + \|y\|^4}{2}
	\end{equation}
	where $\frac{(\lambda^2 + \|x\|^2)^2 + (\lambda^2 + \|y\|^2)^2}{2}$ and $\frac{2\lambda^4 + \|x\|^4 + \|y\|^4}{2}$ are both convex functions. Note that this dc decomposition is different from the one presented in \cite{Niu15}.
	\item In order to get a dc decomposition for $f_1$, we can rewrite $\lambda$ as the dc decomposition
	\begin{equation}
	\lambda = \frac{(\lambda + 1)^2}{2} - \frac{\lambda^2+1}{2}
	\end{equation}
	and $y^{\top}(x+z)$ as the dc decomposition
	\begin{equation}
	y^{\top}(x+z) = \biggl(\frac{\|y+x\|^2}{4} + \frac{\|y+z\|^2}{4}\biggr) -  \biggl(\frac{\|y-x\|^2}{4} + \frac{\|y-z\|^2}{4}\biggr).
	\end{equation}
	Finally, $f_1$ is the product of two dc functions $-2\lambda$ and $y^{\top}(x+z)$ with positive dc components. The following proposition gives an explicit dc decomposition of the product function $f_1f_2$.
	\begin{prop}[\cite{Horst}]\label{prop:dcdecompofproduct}
		If $\phi_1=g_1-h_1$ and $\phi_2=g_2-h_2$ whose dc components $g_1,h_1,g_2,h_2$ are all positive convex functions, then $\phi_1\phi_2$ has a dc decomposition
		\[\phi_1\phi_2 = \frac{(g_1+g_2)^2 + (h_1+h_2)^2}{2} - \frac{(g_1+h_2)^2 + (g_2+h_1)^2}{2}. \] 
	\end{prop}
	Based on Proposition \ref{prop:dcdecompofproduct}, we obtain a dc decomposition of $f_1$ as:
	\begin{multline}\label{eq:dcdf1}
	f_1(x,y,z,\lambda) = \\
	\frac{(4\lambda^2+4 + \|y+x\|^2 + \|y+z\|^2)^2 + (4(\lambda+1)^2 + \|y-x\|^2 + \|y-z\|^2)^2}{32}\\
	- \frac{(4\lambda^2+4 + \|y-x\|^2 + \|y-z\|^2)^2 + (4(\lambda+1)^2 + \|y+x\|^2 + \|y+z\|^2)^2}{32}.
	\end{multline}
\end{enumerate}

From (\ref{eq:df})--(\ref{eq:dcdf2}) and (\ref{eq:dcdf1}), a dc decomposition for the polynomial objective function $f$ is given by:
\begin{equation}
f(x,y,z,w,\lambda) = g(x,y,z,w,\lambda) - h(x,y,z,w,\lambda)
\end{equation}
where 
\begin{equation}\label{eq:g&h1}
\left\{\begin{array}{ll}
g(x,y,z,w,\lambda) = &\|y\|^2 + \|z\|^2 + \frac{\|x+w\|^2}{4} + \frac{(\lambda^2 + \|x\|^2)^2 + (\lambda^2 + \|y\|^2)^2}{2} \\
& + \frac{(4\lambda^2+4 + \|y+x\|^2 + \|y+z\|^2)^2 + (4(\lambda+1)^2 + \|y-x\|^2 + \|y-z\|^2)^2}{32},\\
h(x,y,z,w,\lambda) = &\frac{\|x-w\|^2}{4} + \frac{2\lambda^4 + \|x\|^4 + \|y\|^4}{2} \\
&+ \frac{(4\lambda^2+4 + \|y-x\|^2 + \|y-z\|^2)^2 + (4(\lambda+1)^2 + \|y+x\|^2 + \|y+z\|^2)^2}{32}.
\end{array}
\right.
\end{equation}

Hence, a dc programming formulation of $(P)$ is finally stated as follows:

\begin{Dcp}{($P_{DC}$)}{pdc}
	\begin{equation}\label{prob:pdc}
	 \min \{g(x,y,z,w,\lambda)-h(x,y,z,w,\lambda): (x,y,z,w,\lambda)\in \mathcal{C}\}
	\end{equation} 
	where $g$ and $h$ are defined in (\ref{eq:g&h1}) and $\mathcal{C}$ is given by (\ref{eq:C}).
\end{Dcp}

Note that the dc components $g$ and $h$ defined by (\ref{eq:g&h1}) are written in form of sums-of-squares (SOS). This decomposition is called \emph{Difference of Convex SOS} (DCSOS) decomposition, which is an extension of the Difference of SOS (DSOS) decomposition introduced by Niu in \cite{Niu17}. It is proved that any polynomial can be rewritten as DCSOS whose decomposition can be constructed in polynomial time. Another dc decomposition for polynomial functions is the Difference of SOS-convex decomposition proposed by Ahmadi and Hall \cite{Ahmadi} whose dc components are not supposed to be positive. The reader can refer to \cite{Niu17,Ahmadi} for more discussions about DCSOS and Difference of SOS-convex decompositions for polynomials.

\section{A second dc programming formulation}\label{sec:dcpofPhat}

This dc decomposition is based on an \emph{universal dc decomposition technique} for the functions $f_1$ and $f_2$. This approach exploits the fact that the spectral radius of the Hessian matrices of $f_1$ and $f_2$ are bounded on the convex set $\mathcal{C}$. In fact, if $\phi: D\to \R$ is a convex function of class $\mathcal{C}^2$ defined on a compact convex set $D$ of $\R^n$, then the spectral radius of the Hessian matrix $\nabla^2 \phi(x)$, denoted by $\rho(\nabla^2 \phi(x))$, is bounded on $D$. We get a dc decomposition for $\phi$ as:
\begin{equation}\label{eq:udcd}
\phi(x) = \biggl(\frac{\rho^*}{2}\|x\|^2\biggr) - \biggl(\frac{\rho^*}{2}\|x\|^2 - \phi(x)\biggr)
\end{equation}  
where $\rho^*$ is a constant verifying $\rho^* \geq \displaystyle\max_{x\in D}\rho(\nabla^2 \phi(x))$. Obviously, $\frac{\rho^*}{2}\|x\|^2$ and $\frac{\rho^*}{2}\|x\|^2 - \phi(x)$ are both convex functions on $D$.

As stated in \cite{Niu15,Niu11}, in order to obtain a dc decomposition for the nonconvex homogeneous polynomial functions $f_1$ and $f_2$, we first compute their Hessian matrices as follows:
\begin{enumerate}
	\item Gradient and Hessian of $f_1$:
	\[
	\nabla f_1(x,y,z,\lambda) = \begin{bmatrix}
	\nabla_x f_1(x,y,z,\lambda)\\
	\nabla_y f_1(x,y,z,\lambda)\\
	\nabla_z f_1(x,y,z,\lambda)\\
	\nabla_\lambda f_1(x,y,z,\lambda)
	\end{bmatrix} = \begin{bmatrix}
	-2\lambda y  \\
	-2\lambda (x+z)\\
	-2\lambda y\\
	-2y^{\top}(x+z)
	\end{bmatrix},
	\]
	
	\[\nabla^2 f_1(x,y,z,\lambda) = \begin{bmatrix}
	0 & -2\lambda I & 0 & -2y \\
	-2\lambda I & 0 & -2\lambda I & -2(x+z) \\
	0 & -2\lambda I & 0 & -2y \\
	-2y^{\top} & -2(x+z)^{\top} & -2y^{\top} & 0
	\end{bmatrix}.
	\]
	\item Gradient and Hessian of $f_2$:
	\[
	\nabla f_2(x,y,\lambda) = \begin{bmatrix}
	\nabla_x f_2(x,y,\lambda)\\
	\nabla_y f_2(x,y,\lambda)\\
	\nabla_\lambda f_2(x,y,\lambda)
	\end{bmatrix} = \begin{bmatrix}
	2\lambda^2 x\\
	2\lambda^2 y\\
	2\lambda(\|x\|^2+\|y\|^2)
	\end{bmatrix},
	\]
	\[\nabla^2 f_2(x,y,\lambda) = \begin{bmatrix}
	2\lambda^2 I & 0 & 4\lambda x \\
	0 & 2\lambda^2 I & 4\lambda y \\
	4\lambda x^{\top} & 4\lambda y^{\top} & 2(\|x\|^2+\|y\|^2)
	\end{bmatrix}.
	\]	
\end{enumerate}

Let us denote $\|x\|_1$ and $\|x\|_{\infty}$ the 1-norm and infinity norm of the vector $x\in \R^n$ respectively. Let $\rho(\nabla^2 f_1):\R^{3n+1}\to \R$ and $\rho(\nabla^2 f_2):\R^{2n+1}\to \R$ be functions defined by
$$\rho(\nabla^2 f_1)(x,y,z,\lambda) = \rho(\nabla^2 f_1(x,y,z,\lambda)),$$ $$\rho(\nabla^2 f_2)(x,y,\lambda) = \rho(\nabla^2 f_2(x,y,\lambda)).$$
The bounds of $\rho (\nabla^2 f_1)$ and $\rho (\nabla^2 f_2)$ on $\mathcal{C}$ are given by the following theorem:
\begin{thm} \label{thm:specoff1&f2}
	$\rho (\nabla^2 f_1)$ and $\rho (\nabla^2 f_2)$ are bounded on $\mathcal{C}$ by 
	$$\rho (\nabla^2 f_1) \leq 2 \max\{|\lambda| + \|y\|_{\infty}, \|x+z\|_{\infty}+2|\lambda|, 2\|y\|_1 + \|x+z\|_{1}\},$$
	$$\rho (\nabla^2 f_2) \leq 2 \max\{\lambda^2+2\|\lambda x\|_{\infty}, \lambda^2+2\|\lambda y\|_{\infty}, \|(x,y)\|^2 + 2(\|\lambda x\|_1 + \|\lambda y\|_1) \}.$$
\end{thm}
\begin{pf}
	These inequalities are due to the fact that $\rho(A)\leq \|A\|_1$ where $\|A\|_1$ stands for the induced $1$-norm of the matrix $A$ defined by $\|A\|_1 = \max_{j\in \IntEnt{1,n}}\sum_{i=1}^{n}|A_{ij}|$ (in which $\IntEnt{1,n}$ stands for the set $\{1,\ldots, n\}$). Thus it is sufficient to compute $\|\nabla^2 f_1\|_1$ and $\|\nabla^2 f_2\|_1$ with $(x,y,z,\lambda)\in\mathcal{C}$, which yield the required inequalities.
\end{pf}

It follows from Theorem \ref{thm:specoff1&f2} that the functions $\rho(\nabla^2 f_1)$ and $\rho(\nabla^2 f_2)$ are bounded on a bounded set $\mathcal{C}$. It is known \cite{Fernandes14} that the number of complementary eigenvalues of QEiCP($A,B,C$) is finite. So, if QEiCP($A,B,C$) has a solution, then there exist $l$ and $u$ in $\R$ such that $$l\leq \lambda\leq u$$
for all complementary eigenvalues $\lambda$.
\begin{thm}\label{thm:boundsofvars}
	If QEiCP($A,B,C$) has a solution $(\lambda, x)$ with $\lambda\in [l,u]$, then any optimal solution of $(P)$ satisfies:
	\[x\in [0,1]^n, y\in [\min\{0,l\},\max\{0,u\}]^n, z\in [0,p^2]^n, e^{\top}z\leq p^2, w\in [0,\bar{w}]\]
	with $p = \max\{|l|,|u|\}$ and $\bar{w}=[\bar{w_i}]_{i\in \IntEnt{1,n}}$ is given by \[\bar{w}_i=p^2\sqrt{\sum_{j=1}^{n}A_{ij}^2} + p\sqrt{\sum_{j=1}^{n}B_{ij}^2} +  \sqrt{\sum_{j=1}^{n}C_{ij}^2}, \forall i \in \IntEnt{1,n}.
	\]
\end{thm}
\begin{pf} The proofs for the bounds of $x,y,z$ are similar to the proposition $1$ in \cite{Niu15}. So we only have to prove that $e^{\top}z\leq p^2$ and the new upper bound $\bar{w}$ for $w$. 
	\begin{enumerate}
		\item[(i)] Based on $z=\lambda^2x, e^{\top}x=1$ and $\lambda \in [l,u]$, it follows that $e^{\top}z=\lambda^2\leq p^2$.
		\item[(ii)] The upper bound of $w$ can be obtained from the definition of $w$ as $Az+By+Cx$. Since $x,y,z$ are all bounded, then $w_i, \forall i\in \IntEnt{1,n}$ is also bounded. In fact,
		\begin{align}
		w_i &\leq |w_i| = \left|\sum_{j=1}^{n}A_{ij}z_j + \sum_{j=1}^{n}B_{ij}y_j + \sum_{j=1}^{n}C_{ij}x_j \right|\nonumber\\
		&\leq \left|\sum_{j=1}^{n}A_{ij}z_j\right| + \left|\sum_{j=1}^{n}B_{ij}y_j\right| + \left|\sum_{j=1}^{n}C_{ij}x_j \right|.\nonumber
		\end{align}
		According to the Cauchy-Schwartz inequality, we get respectively that
		\[\left|\sum_{j=1}^{n}A_{ij}z_j\right| \leq \|z\|\sqrt{\sum_{j=1}^{n}A_{ij}^2}=\lambda^2\|x\|\sqrt{\sum_{j=1}^{n}A_{ij}^2}\leq \max\{l^2,u^2\}\sqrt{\sum_{j=1}^{n}A_{ij}^2},\]
		\[\left|\sum_{j=1}^{n}B_{ij}y_j\right| \leq \|y\|\sqrt{\sum_{j=1}^{n}B_{ij}^2}=\|\lambda|\|x\|\sqrt{\sum_{j=1}^{n}B_{ij}^2}\leq \max\{|l|,|u|\}\sqrt{\sum_{j=1}^{n}B_{ij}^2},\]
		and 
		\[\left|\sum_{j=1}^{n}C_{ij}x_j\right| \leq \|x\|\sqrt{\sum_{j=1}^{n}C_{ij}^2}\leq \sqrt{\sum_{j=1}^{n}C_{ij}^2}.\]
		Thus, we have $\forall i\in \IntEnt{1,n}$, 
		\begin{align}
		w_i &\leq \max\{l^2,u^2\}\sqrt{\sum_{j=1}^{n}A_{ij}^2} + \max\{|l|,|u|\}\sqrt{\sum_{j=1}^{n}B_{ij}^2} +  \sqrt{\sum_{j=1}^{n}C_{ij}^2}\nonumber\\
		&=p^2\sqrt{\sum_{j=1}^{n}A_{ij}^2} + p\sqrt{\sum_{j=1}^{n}B_{ij}^2} +  \sqrt{\sum_{j=1}^{n}C_{ij}^2} = \bar{w_i}.\nonumber
		\end{align}
	\end{enumerate}
\end{pf}
\begin{rmk}
	The vector $\bar{w}$ provides a tighter bound for $w$ than the one given in \cite{Niu15}, and its computation requires the knowledge of the bounds $l$ and $u$ for $\lambda$. Formulas for computing these two bounds $l$ and $u$ are introduced in section \ref{sec:boundoflambda}.
\end{rmk}

Since the solution set of $(P)$ is bounded, if we add the bound obtained in Theorem \ref{thm:boundsofvars} into the problem $(P)$, the resulting problem should be equivalent to $(P)$. Let us consider the compact convex polyhedral set 
\begin{align}\label{eq:Chat}
\hat{\mathcal{C}} &= \{(x,y,z,w,\lambda): w=Az+By+Cx, e^{\top}x=1, e^{\top}y=\lambda, e^{\top}z\leq p^2,\lambda \in [l,u], \\
\nonumber & x\in[0,1]^n, w\in [0,\bar{w}], y\in [\min\{0,l\},\max\{0,u\}]^n, z\in [0,p^2]^n\}. 
\end{align}
So the following problem 
\begin{equation}\label{prob:nlphat}
(\hat{P}) \qquad 
\min \{f(x,y,z,w,\lambda) : (x,y,z,w,\lambda) \in \mathcal{\hat{C}}.\}\nonumber
\end{equation}
has the same optimal solution set as $(P)$ when $\rho (\nabla^2 f_1)$ and $\rho (\nabla^2 f_2)$ are bounded.

If we assume that $C\notin S_0$, then $\lambda=0$ cannot be a solution of QEiCP and this problem has at least a positive and a negative complementarity eigenvalues \cite{Iusem16}. So $\hat{C}$ can be reduced into two subsets regarding to $\lambda$ positive or negative. If $\lambda\geq 0$, then $y=\lambda x\geq 0$, and $(x,y,z,w,\lambda)$ in $\hat{C}$ should be nonnegative. Otherwise, $\lambda\leq 0$, then $y\leq 0$ and $(x,z,w)$ is still nonnegative. Let us define the two subsets
\begin{align}\label{eq:Chat1}
\hat{\mathcal{C}}_1 &= \{(x,y,z,w,\lambda): w=Az+By+Cx, e^{\top}x=1, e^{\top}y=\lambda, \\\nonumber 
& e^{\top}z\leq p^2,\lambda \in [0,u], x\in[0,1]^n, w\in [0,\bar{w}], y\in [0,u]^n, z\in [0,p^2]^n\}. 
\end{align}
\begin{align}\label{eq:Chat2}
\hat{\mathcal{C}}_2 &= \{(x,y,z,w,\lambda): w=Az+By+Cx, e^{\top}x=1, e^{\top}y=\lambda, \\
\nonumber &  e^{\top}z\leq p^2,\lambda \in [l,0], x\in[0,1]^n, w\in [0,\bar{w}], y\in [l,0]^n, z\in [0,p^2]^n\}. 
\end{align}

Obviously, problem $(\hat{P})$ is equivalent to problem:
\begin{equation}
\min \{f(x,y,z,w,\lambda) : (x,y,z,w,\lambda) \in \mathcal{\hat{C}}_1 \cup \mathcal{\hat{C}}_2 \}
\end{equation}
which means that, for seeking a positive (resp. negative) complementary eigenvalue, we just need to solve the problem $(\hat{P})$ on $\hat{\mathcal{C}}_1$ (resp. $\hat{\mathcal{C}}_2$).

The next theorem provides upper bounds for $\rho (\nabla^2 f_1)$ and $\rho (\nabla^2 f_2)$ on $\mathcal{\hat{C}}_1$ and $\mathcal{\hat{C}}_2$.

\begin{thm}\label{thm:boundsofspectralradius} On $\mathcal{\hat{C}}_1$ and $\mathcal{\hat{C}}_2$, we have
	\[\rho (\nabla^2 f_1) \leq \rho_1, \rho (\nabla^2 f_2) \leq \rho_2\]
	with $\rho_1= 2(p+1)^2$, $\rho_2 = 6p^2+4p+2$ and $p=\max\{|l|,|u|\}$. 
\end{thm}
\begin{pf}
	Since $|\lambda| \leq  \max\{|l|,|u|\}=p$ both on $\mathcal{\hat{C}}_1$ and $\mathcal{\hat{C}}_2$, it follows from Theorem \ref{thm:specoff1&f2} that
	$\rho (\nabla^2 f_1) \leq  2 \max\{|\lambda| + \|y\|_{\infty}, \|x+z\|_{\infty}+2|\lambda|, 2\|y\|_1 + \|x+z\|_{1}\}.$
	According to Theorem \ref{thm:boundsofvars}, we get
	\[|\lambda|\leq p, \|x\|_1 =1, \|x\|\leq 1, \|x\|_{\infty}\leq 1.\]
	\[\|y\|_1 =\sum_{i=1}^{n}|y_i| = \biggl(\begin{cases}
	e^{\top}y & \text{, on } \mathcal{\hat{C}}_1\\
	-e^{\top}y & \text{, on } \mathcal{\hat{C}}_2
	\end{cases}\biggr) \leq |\lambda| \leq p. \]
	\[\|y\| \leq p, \|y\|_{\infty}\leq p, \|z\|_1\leq p^2, \|z\|_{\infty}\leq p^2.\]
	\[\|x+z\|_{\infty} \leq \|x\|_{\infty} + \|z\|_{\infty} \leq 1 + p^2 \text{ and } \|x+z\|_{1} \leq \|x\|_{1} + \|z\|_{1} \leq 1 + p^2.\]	
	Hence,
	\[\rho (\nabla^2 f_1)\leq 2\max\{2p,1+p^2+2p,1+p^2+2p\}= 2(p+1)^2=\rho_1.\]
	Similarly, 
	\[\rho (\nabla^2 f_2) \leq 2 \max\{\lambda^2+2\|\lambda x\|_{\infty}, \lambda^2+2\|\lambda y\|_{\infty}, \|x\|^2+\|y\|^2 + 2(\|\lambda x\|_1 + \|\lambda y\|_1) \} \]
	\[\leq 2\max\{p^2+2p, 3p^2, 3p^2+2p+1\} =6p^2+4p+2=\rho_2.\]
\end{pf}
\begin{rmk}
	These upper bounds $\rho_1$ and $\rho_2$ for spectral radius are tighter than those given in our previous work \cite{Niu15}.
\end{rmk}

Finally, we can obtain a similar dc decomposition for $f_1$ and $f_2$ as in \cite{Niu15} but with smaller values for $\rho_1$ and $\rho_2$ as:
\begin{equation}\label{eq:dcdf1withspectralradius}
f_1(x,y,z,\lambda)=\frac{\rho_1}{2} \|(x,y,z,\lambda)\|^2 - \left(\frac{\rho_1}{2} \|(x,y,z,\lambda)\|^2 - f_1(x,y,z,\lambda)\right).
\end{equation}
\begin{equation}\label{eq:dcdf2withspectralradius}
f_2(x,y,\lambda)=\frac{\rho_2}{2} \|(x,y,\lambda)\|^2 - \left(\frac{\rho_2}{2} \|(x,y,\lambda)\|^2 - f_2(x,y,\lambda)\right).
\end{equation}
where $ \frac{\rho_1}{2} \|(x,y,z,\lambda)\|^2 - f_1(x,y,z,\lambda)$ and $\frac{\rho_2}{2} \|(x,y,\lambda)\|^2 - f_2(x,y,\lambda)$ are both locally convex functions on $\mathcal{\hat{C}}_1$ and $\mathcal{\hat{C}}_2$, and $\frac{\rho_1}{2} \|(x,y,z,\lambda)\|^2$ and $\frac{\rho_2}{2} \|(x,y,\lambda)\|^2$ are quadratic convex functions.

Hence, we get from  (\ref{eq:df})--(\ref{eq:dcdf3}), (\ref{eq:dcdf1withspectralradius}) and (\ref{eq:dcdf2withspectralradius}) a dc decomposition of $f$ as:
\[f(x,y,z,w,\lambda) = \hat{g}(x,y,z,w,\lambda)-\hat{h}(x,y,z,w,\lambda)\]
with
\begin{equation}\label{eq:ghatandhhat}
\left\lbrace \begin{array}{l}
\hat{g}(x,y,z,w,\lambda) = \frac{\|x+w\|^2}{4} + \|(y,z)\|^2+ \frac{\rho_1}{2} \|(x,y,z,\lambda)\|^2  + \frac{\rho_2}{2} \|(x,y,\lambda)\|^2, \\
\hat{h}(x,y,z,w,\lambda) = \hat{g}(x,y,z,w,\lambda)-f(x,y,z,w,\lambda),
\end{array}\right.
\end{equation}
which yields the following dc programming formulation of $(\hat{P})$:

\begin{Dcp}{($\hat{P}_{DC}$)}{pdchat}
	\begin{equation}\label{prob:pdchat}
	 \min \{\hat{g}(x,y,z,w,\lambda)-\hat{h}(x,y,z,w,\lambda): (x,y,z,w,\lambda)\in \mathcal{\hat{C}}_1 (\text{or } \mathcal{\hat{C}}_2)\}
	\end{equation} 
	with $\hat{g}$ and $\hat{h}$ defined in (\ref{eq:ghatandhhat}), and $\mathcal{\hat{C}}_1$ and $\mathcal{\hat{C}}_2$ are given by (\ref{eq:Chat1}) and (\ref{eq:Chat2}).
\end{Dcp}

\section{DC programming formulations of $(P')$}\label{sec:dcpofP'}
In a similar fashion to the two previous sections, we can obtain two different dc programming formulations for $(P')$. The only difference between $(P)$ and $(P')$ is related with the last term of the objective function where the bilinear term $f_3(x,w) = w^{\top}x$ is replaced by the concave polyhedral function $\tilde{f}_3(x,w) = \sum_{i=1}^{n} \min\{x_i,w_i\}$. This concave function $\tilde{f_3}$ has an explicit dc decomposition as: 
\begin{equation}\label{eq_f3tilde}
\tilde{f_3}(x,w) = \sum_{i=1}^{n} \min\{x_i,w_i\}= (0) - (-\sum_{i=1}^{n} \min\{x_i,w_i\})
\end{equation}
where $-\sum_{i=1}^{n} \min\{x_i,w_i\}$ is a convex polyhedral function.

Thus we can do exactly in the same way as in $(P)$ and $(\hat{P})$ by only replacing $f_3$ by $\tilde{f}_3$ to get the nonlinear programming formulations $(P')$ and $(\hat{P}')$ as well as their corresponding dc programming formulations.

For problem
\begin{equation}
(P') \qquad 
\min \{ f'(x,y,z,w,\lambda) : (x,y,z,w,\lambda)\in \mathcal{C}\},\nonumber
\end{equation}
we can generate a dc decomposition for $f'$ as 
\begin{equation}
f'(x,y,z,w,\lambda) = g'(x,y,z,w,\lambda) - h'(x,y,z,w,\lambda)
\end{equation}
where 
\begin{equation}
\left\{\begin{array}{ll}\label{eq:g'&h'}
g'(x,y,z,w,\lambda) = &\|y\|^2 + \|z\|^2 +  \frac{(\lambda^2 + \|x\|^2)^2 + (\lambda^2 + \|y\|^2)^2}{2} \\
& + \frac{(4\lambda^2+4 + \|y+x\|^2 + \|y+z\|^2)^2 + (4(\lambda+1)^2 + \|y-x\|^2 + \|y-z\|^2)^2}{32},\\
h'(x,y,z,w,\lambda) = &-\sum_{i=1}^{n} \min\{x_i,w_i\} + \frac{2\lambda^4 + \|x\|^4 + \|y\|^4}{2} \\
&+ \frac{(4\lambda^2+4 + \|y-x\|^2 + \|y-z\|^2)^2 + (4(\lambda+1)^2 + \|y+x\|^2 + \|y+z\|^2)^2}{32}.
\end{array}
\right.
\end{equation}
Hence, a dc programming formulation of $(P')$ is given by:

\begin{Dcp}{($P'_{DC}$)}{pdc'}
	\begin{equation}\label{prob:pdc'}
	 \min \{g'(x,y,z,w,\lambda)-h'(x,y,z,w,\lambda): (x,y,z,w,\lambda)\in \mathcal{C}\}
	\end{equation} 
	where $\mathcal{C}$ is defined in (\ref{eq:C}) and $g'$ and $h'$ are given by (\ref{eq:g'&h'}).
\end{Dcp}

We can also get an equivalent problem $(\hat{P}')$ similar to $(\hat{P})$ as:
\begin{equation}\label{prob:nlphat'}
(\hat{P}') \qquad 
\min \{ f'(x,y,z,w,\lambda) : (x,y,z,w,\lambda) \in \mathcal{\hat{C}}\}.\nonumber
\end{equation}

By the universal dc decomposition, a dc decomposition for $f'$ is given by:
\[f'(x,y,z,w,\lambda) = \hat{g}'(x,y,z,w,\lambda)-\hat{h}'(x,y,z,w,\lambda)\]
with
\begin{equation}\label{eq:ghat'andhhat'}
\left\lbrace \begin{array}{l}
\hat{g}'(x,y,z,w,\lambda) = \|(y,z)\|^2+ \frac{\rho_1}{2} \|(x,y,z,\lambda)\|^2  + \frac{\rho_2}{2} \|(x,y,\lambda)\|^2, \\
\hat{h}'(x,y,z,w,\lambda) = \hat{g}'(x,y,z,w,\lambda)-f'(x,y,z,w,\lambda),
\end{array}\right. 
\end{equation}
which yields the following dc programming formulation of $(\hat{P}')$:

\begin{Dcp}{($\hat{P}'_{DC}$)}{pdchat'}
	\begin{equation}\label{prob:pdchat'}
	 \min \{\hat{g}'(x,y,z,w,\lambda)-\hat{h}'(x,y,z,w,\lambda): (x,y,z,w,\lambda)\in \mathcal{\hat{C}}_1 (\text{or } \mathcal{\hat{C}}_2)\}
	\end{equation} 
	with $\hat{g}'$ and $\hat{h}'$ defined in (\ref{eq:ghat'andhhat'}), and $\mathcal{\hat{C}}_1$ and $\mathcal{\hat{C}}_2$ are given by (\ref{eq:Chat1}) and (\ref{eq:Chat2}).
\end{Dcp}

\section{Lower and upper bounds for $\lambda$}\label{sec:boundoflambda}
The dc formulations described in the previous sections require the existence of lower and upper bounds for the variable $\lambda$. In this section, we discuss the computation of such bounds for QEiCP with a matrix $A\in$ PD and satisfying the co-hyperbolic condition or $C\notin S_0$. 

\subsection{Lower and upper bounds of $\lambda$ when $A\in$ PD and co-hyperbolicity}

The bounds are given in the following theorems:

\begin{thm}[see \cite{Niu15}]\label{thm:boundsoflambda}
	If $A\in$ PD and the co-hyperbolic condition holds, the $\lambda$-component of any solution of QEiCP satisfies
	\[l=\beta-\sqrt{\alpha} \leq \lambda \leq \gamma+\sqrt{\alpha}=u\]
	with $$\alpha = \max\{\gamma^2,\beta^2\} + \max_{i,j}\{-C_{ij}\}/s,$$
	\[\beta = \left\{\begin{array}{ll}
	\min\{-B_{ij}\}/(2\max\{A_{ij}\}), & \hbox{if $\min\{-B_{ij}\}>0$;} \\
	\min\{-B_{ij}\}/(2s), & \hbox{if $\min\{-B_{ij}\}\leq 0$,}
	\end{array}
	\right.\]
	\[\gamma=\left\{\begin{array}{ll}
	\max\{-B_{ij}\}/(2s), & \hbox{if $\max\{-B_{ij}\}>0$;} \\
	\max\{-B_{ij}\}/(2\max\{A_{ij}\}), & \hbox{if $\max\{-B_{ij}\}\leq0$,}
	\end{array}
	\right.\]
	in which $s=\min\{x^{\top}Ax: e^{\top}x=1, x\geq 0\}$. 
\end{thm}
\begin{pf}
	The reader is referred to the proof given in \cite{Niu15}.
\end{pf}

Based on a similar idea used in Theorem \ref{thm:boundsoflambda}, we get a new bound for $\lambda$ in the next result:
\begin{thm}\label{thm:boundsoflambdabis}
	Suppose that $A\in$ PD and QEiCP satisfies the co-hyperbolic condition. Let $\lambda_{\max}(M)$ (resp. $\lambda_{\min}(M)$) be the largest (resp. the smallest) real eigenvalue of the matrix $M$. The $\lambda$-component of any solution of QEiCP satisfies
	\[l=\beta-\sqrt{\alpha} \leq \lambda \leq \gamma+\sqrt{\alpha}=u\]
	with $a=\lambda_{\min}(A+A^{\top})/2$, $\bar{a}=\lambda_{\max}(A+A^{\top})/2$, $b=\lambda_{\min}(-B-B^{\top})/2$,  $\bar{b}=\lambda_{\max}(-B-B^{\top})/2$, $\bar{c}=\lambda_{\max}(-C-C^{\top})/2$, $\alpha = \max\{\gamma^2,\beta^2\} + \bar{c}/a$,
	\[\beta = \left\{\begin{array}{ll}
	b/(2\bar{a}) &, \hbox{if $b>0$;} \\
	b/(2a) &, \hbox{if $b\leq 0$,}
	\end{array}
	\right.
	\gamma=\left\{\begin{array}{ll}
	b/(2a) &, \hbox{if $\bar{b}>0$;} \\
	b/(2\bar{a}) &, \hbox{if $\bar{b}\leq0$.}
	\end{array}
	\right.\]
\end{thm}
\begin{pf}
	The proof is almost the same as in Theorem \ref{thm:boundsoflambda}. The major difference is to use the set $V=\{x\in \R^n: \|x\|=1,x\geq 0 \}$ instead of the set $U$ used in the proof of this theorem \cite{Niu15}. The reason for introducing the constraint $e^{\top}x=1$ is to avoid a null vector to be a solution of QEiCP. Obviously, $\forall M\in \M_n$ and $x\in V$, we have
	\[x^{\top}Mx = x^{\top} \frac{M+M^{\top}}{2}x.\]
	Since $(M+M^{\top})/2$ is a symmetric matrix, then it is diagonalizable and the quadratic form defined on $V$ is bounded by 
	\[\lambda_{\min}\big( \frac{M+M^{\top}}{2}\big) \leq x^{\top} \frac{M+M^{\top}}{2}x \leq \lambda_{\max}\big( \frac{M+M^{\top}}{2}\big), \forall x\in V.  \]
	Therefore, we can use the above inequality to get our new bounds of $\lambda$ following the same procedure used in the proof of Theorem \ref{thm:boundsoflambda}.
\end{pf}
\begin{rmk}
	Note that another bound of $\lambda$ has been provided in \cite{Fernandes14}, and it seems interesting to compare these two bounds. Theorem \ref{thm:boundsofvars} shows that the bounds of $(x,y,z,w)$ depend on the bounds of $\lambda$, and the tighter bound of $\lambda$ is, the tighter the bounds of $(w,y,z,w)$ are.
\end{rmk}

\subsection{Lower and upper bounds of $\lambda$ for $A\in$ PD and $C\notin S_0$}\label{subsec:CnotinS0}
In this case, it is proved that QEiCP($A,B,C$) has both positive and negative complementary eigenvalues \cite{Bras16}. We just need to focus on the bounds for a positive complementary eigenvalue, since the negative one can be done in a similar way. 
\begin{thm}[\cite{Iusem16}]
	Let $A\in$ PD, $C\notin S_0$ and define the following programs
	$$(LP) \qquad \min\{e^{\top}v + e^{\top}y : Av + By + Cx \geq 0, e^{\top}y + e^{\top}x = 1, (x,y,v) \geq 0\}$$ 
	and
	$$(UP) \qquad \max\{ \frac{p^{\top}y}{y^{\top}Ay + x^{\top}x} : e^{\top}y + e^{\top}x = 1, (x,y)\geq 0\}$$
	where $p$ is a vector with components $p_i = 1+\sum_{j=1}^{n}(\max\{0,−B_{ij}\}+\max\{0,−C_{ij}\}), \forall i\in \IntEnt{1,n}$. Then  
	\begin{itemize}
		\item[(i)] The linear program $(LP)$ has an optimal solution with positive optimal value which is a lower bound of $\lambda$.
		\item[(ii)] Any stationary point of $(UP)$ is a global maximum of $(UP)$ and its optimal value is an upper bound of $\lambda$.
	\end{itemize}
\end{thm}

\section{DC algorithms for solving dc programming formulations}\label{sec:dcas}
In this section, we investigate how to use DCA for solving the dc programming formulations $(P_{DC})$, $(\hat{P}_{DC})$, $(P'_{DC})$ and $(\hat{P}'_{DC})$.  
\subsection{DCA for solving $(P_{DC})$}
Since the function $h$ defined in (\ref{eq:g&h1}) is differentiable, then $\partial h(x,y,z,w,\lambda)$ is reduced to a singleton $\{\nabla h(x,y,z,w,\lambda)\}$. The following fixed-point scheme describes the major computations in dc algorithm for $(P_{DC})$ \begin{equation}\label{prob:dcaforpdc}
\begin{array}{l}
(x^{k+1},y^{k+1},z^{k+1},w^{k+1},\lambda^{k+1}) \in \argmin\{g(x,y,z,w,\lambda) \\
-\langle (x,y,z,w,\lambda),\nabla h(x^k,y^k,z^k,w^k,\lambda^k)\rangle: (x,y,z,w,\lambda)\in \mathcal{C}\}
\end{array}
\end{equation}
where $g$ is defined in (\ref{eq:g&h1}), $\mathcal{C}$ is defined in (\ref{eq:C}), and $\langle u,v \rangle$ represents inner product of two vectors $u$ and $v$. Since this problem is a convex polynomial optimization, then any KKT solution is a global optimal solution. 

The gradient of $h$ can be computed as follows:
\begin{equation}\label{eq:dh}
\nabla h(x,y,z,w,\lambda)  = \begin{bmatrix}
\nabla_x h(x,y,z,w,\lambda)\\
\nabla_y h(x,y,z,w,\lambda)\\
\nabla_z h(x,y,z,w,\lambda)\\
\nabla_w h(x,y,z,w,\lambda)\\
\nabla_{\lambda} h(x,y,z,w,\lambda)\\
\end{bmatrix}
\end{equation}
where
\begin{equation}\label{eq:dhforPdc}
\left\lbrace 
\begin{array}{ll}
\nabla_x h(x,y,z,w,\lambda) = &\frac{x-w}{2}+ 2\|x\|^2x + \frac{(4\lambda^2+4+\|y-x\|^2+\|y-z\|^2)(x-y)}{8} \\
&+ \frac{(4(\lambda+1)^2+\|y+x\|^2+\|y+z\|^2)(x+y)}{8}.\\
\nabla_y h(x,y,z,w,\lambda) = &2\|y\|^2y + \frac{(4\lambda^2+4+\|y-x\|^2+\|y-z\|^2)(2y-x-z)}{8} \\
&+ \frac{(4(\lambda+1)^2+\|y+x\|^2+\|y+z\|^2)(2y+x+z)}{8}.\\
\nabla_z h(x,y,z,w,\lambda) = & \frac{(4\lambda^2+4+\|y-x\|^2+\|y-z\|^2)(z-y)}{8} \\
&+ \frac{(4(\lambda+1)^2+\|y+x\|^2+\|y+z\|^2)(z+y)}{8}.\\
\nabla_w h(x,y,z,w,\lambda) = &\frac{w-x}{2}.\\
\nabla_{\lambda} h(x,y,z,w,\lambda) = &4\lambda^3 + \frac{(4\lambda^2+4+\|y-x\|^2+\|y-z\|^2)\lambda}{2} \\
&+ \frac{(4(\lambda+1)^2+\|y+x\|^2+\|y+z\|^2)(\lambda+1)}{2}.
\end{array}
\right. 
\end{equation}

The steps of DCA for solving $(P_{DC})$ are described in Algorithm \ref{alg:dcaforpdc}:
\begin{algorithm}[H]
	\caption{DCA for solving $(P_{DC})$}
	\label{alg:dcaforpdc}
	\begin{algorithmic}
		\STATE{\textbf{Inputs}: Initial point $(x^0,y^0,z^0,w^0,\lambda^0)$, the tolerance for optimal value $\epsilon_1 >0$, the tolerance for optimal solution $\epsilon_2>0$, and the tolerance for globality $\epsilon_3>0$.}
		\STATE{\textbf{Outputs}: Optimal solution $(x^*,y^*,z^*,w^*,\lambda^*)$ and optimal value $f^*$.}
		\STATE{}
		\STATE{Set $k=0$, $f^*= +\infty$, $\Delta f = +\infty$ and $\Delta X = +\infty$.}
		\WHILE{($\Delta f > \epsilon_1$ and $\Delta X > \epsilon_2$ and $f^*>\epsilon_3$)}
		\STATE{
			\begin{enumerate}[\textbf{Step} 1:]
				\item Compute $\nabla h(x^k,y^k,z^k,w^k,\lambda^k)$ by formula (\ref{eq:dhforPdc}).
				\item Compute $(x^{k+1},y^{k+1},z^{k+1},w^{k+1},\lambda^{k+1})$ by solving the convex optimization problem (\ref{prob:dcaforpdc}).
				\item Compute 
				$$\Delta f \leftarrow |f(x^{k+1},y^{k+1},z^{k+1},w^{k+1},\lambda^{k+1}) - f(x^{k},y^{k},z^{k},w^{k},\lambda^{k})|.$$
				$$\Delta X \leftarrow \|(x^{k+1},y^{k+1},z^{k+1},w^{k+1},\lambda^{k+1}) - (x^{k},y^{k},z^{k},w^{k},\lambda^{k})\|.$$
				$$(x^*,y^*,z^*,w^*,\lambda^*)\leftarrow (x^{k+1},y^{k+1},z^{k+1},w^{k+1},\lambda^{k+1}).$$
				$$f^* \leftarrow f(x^{*},y^{*},z^{*},w^{*},\lambda^{*}).$$
				
				Update $k\leftarrow k+1.$
			\end{enumerate}
		}
		\ENDWHILE
		\RETURN $(x^*,y^*,z^*,w^*,\lambda^*)$ and $f^*$.
	\end{algorithmic}
\end{algorithm}

This algorithm should terminate if one of the following stopping criteria is satisfied:
\begin{enumerate}
	\item The sequence $\{f(x^k, y^k ,z^k, w^k,\lambda^k)\}$  converges, i.e., $\Delta f \leq \epsilon_1.$
	\item The sequence $\{(x^k,y^k,z^k,w^k,\lambda^k)\}$ converges, i.e., $\Delta X \leq \epsilon_2.$
	\item The \emph{sufficient global optimality condition} holds, i.e., $f^*\leq \epsilon_3$.
\end{enumerate}

We have the following convergence theorem for DCA:
\begin{thm}[Convergence theorem for DCA]\label{thm:convergenceofdca}
	DCA for solving $(P_{DC})$ generates convergent sequences $\{(x^k,y^k,z^k,w^k,\lambda^k)\}$ and $\{f(x^k,y^k,z^k,w^k,\lambda^k)\}$ such that 
	\begin{itemize}
		\item[(i)] The sequence $\{f(x^k,y^k,z^k,w^k,\lambda^k)\}$ is decreasing and bounded below.
		\item[(ii)] The sequence $\{(x^k,y^k,z^k,w^k,\lambda^k)\}$ converges either to an approximate solution of QEiCP satisfying the sufficient global optimality condition, or to an approximate solution of general KKT point of $(P_{DC})$.
	\end{itemize}
\end{thm}
\begin{pf}
	This theorem is a consequence of the general convergence theorem for DCA \cite{Pham98,Pham05}. The sufficient global optimality condition is based on Theorem \ref{thm:equivofQEiCPandP}, which states that the optimal value of $(P_{DC})$ is zero for any solution of QEiCP.
\end{pf}
\begin{rmk}
	If the sequence $f(x^k,y^k,z^k,w^k,\lambda^k)$ does not converge to zero, it may be caused by one of the following reasons:
	\begin{itemize}
		\item[(i)] QEiCP is infeasible.
		\item[(ii)] The computed solution found by DCA is a KKT point of $(P_{DC})$ which is not a global minimum.
	\end{itemize}
	Note that it is very easy to show whether QEiCP is infeasible since this verification reduces to the solution of a linear program \cite{Cottle}. If QEiCP is feasible and DCA terminates in the second case, then we should combine it with a global optimization solver in order to find a solution of QEiCP. The hybrid enumerative approach method described in \cite{Iusem16} can be modified in order to incorporate with DCA. This should be a topic for future research.
\end{rmk}

\subsection{DCA for solving $(\hat{P}_{DC})$, $(P'_{DC})$ and $(\hat{P}'_{DC})$}
This method is Algorithm \ref{alg:dcaforpdc} with slightly differences in steps $1$ and $2$, which are discussed below.
\paragraph{\textbf{Step 1 for $(\hat{P}_{DC})$}} Since $\hat{h}$ defined in (\ref{eq:ghatandhhat}) is differentiable, $\partial \hat{h}(x,y,z,w,\lambda)$ is reduced to a singleton $\{\nabla \hat{h}(x,y,z,w,\lambda)\}$, which is given by:
\begin{equation}\label{eq:dhhat}
\begin{array}{rl}
\nabla \hat{h}(x,y,z,w,\lambda)  & = \nabla (\hat{g} - f) (x,y,z,w,\lambda)\\
&=\begin{bmatrix}
\frac{x+w}{2}+(\rho_1+\rho_2)x + 2\lambda y - 2\lambda^2 x - w \\
(\rho_1+\rho_2-2\lambda^2)y + 2\lambda(x+z) \\
\rho_1z + 2\lambda y \\
\frac{w-x}{2} \\
(\rho_1+\rho_2-2(\|x\|^2+\|y\|^2))\lambda + 2y^{\top}(x+z)
\end{bmatrix}.
\end{array}
\end{equation}

\paragraph{\textbf{Step 1 for $(P'_{DC})$}} In this case, $h'$ defined in (\ref{eq:g'&h'}) is non-differentiable. Hence $\partial h'(x,y,z,w,\lambda)$ is an nonempty convex set, which is computed by:
\begin{equation}\label{eq:dh'}
\partial h'(x,y,z,w,\lambda)  =\begin{bmatrix}
\partial_x h'(x,y,z,w,\lambda)\\
\partial_y h'(x,y,z,w,\lambda)\\
\partial_z h'(x,y,z,w,\lambda)\\
\partial_w h'(x,y,z,w,\lambda)\\
\partial_{\lambda} h'(x,y,z,w,\lambda)
\end{bmatrix}
\end{equation}
where
\begin{equation}
\left\lbrace 
\begin{array}{ll}
\partial_x h'(x,y,z,w,\lambda) = &-u+ 2\|x\|^2x + \frac{(4\lambda^2+4+\|y-x\|^2+\|y-z\|^2)(x-y)}{8} \\
&+ \frac{(4(\lambda+1)^2+\|y+x\|^2+\|y+z\|^2)(x+y)}{8}.\\
\partial_y h'(x,y,z,w,\lambda) = &2\|y\|^2y + \frac{(4\lambda^2+4+\|y-x\|^2+\|y-z\|^2)(2y-x-z)}{8} \\
&+ \frac{(4(\lambda+1)^2+\|y+x\|^2+\|y+z\|^2)(2y+x+z)}{8}.\\
\partial_z h'(x,y,z,w,\lambda) = & \frac{(4\lambda^2+4+\|y-x\|^2+\|y-z\|^2)(z-y)}{8} \\
&+ \frac{(4(\lambda+1)^2+\|y+x\|^2+\|y+z\|^2)(z+y)}{8}.\\
\partial_w h'(x,y,z,w,\lambda) = &-v.\\
\partial_{\lambda} h'(x,y,z,w,\lambda) = &\frac{(4\lambda^2+4+\|y-x\|^2+\|y-z\|^2)\lambda}{2} \\
&+ \frac{(4(\lambda+1)^2+\|y+x\|^2+\|y+z\|^2)(\lambda+1)}{2}.
\end{array}
\right. 
\end{equation}
and
\begin{equation}\label{eq:u&v}
\left\lbrace \begin{array}{l}
u=[u_i]_{i\in \IntEnt{1,n}}, \text{ with }u_i \in \left\{
\begin{array}{ll}
\{1\}, & \hbox{$x_i< w_i$;} \\
\{0,1\}, &  \hbox{$x_i=w_i$;}\\
\{0\}, & \hbox{$x_i> w_i$.}
\end{array}
\right.\\
v=[v_i]_{i\in \IntEnt{1,n}}, \text{ with } v_i \in \left\{
\begin{array}{ll}
\{0\}, & \hbox{$x_i< w_i$;} \\
\{0,1\}, &  \hbox{$x_i=w_i$;}\\
\{1\}, & \hbox{$x_i> w_i$.}
\end{array}
\right.
\end{array}
\right. 
\end{equation}

\paragraph{\textbf{Step 1 for $(\hat{P}'_{DC})$}}  The function $\hat{h}'$ defined in (\ref{eq:ghat'andhhat'}) is also  non-differentiable. So $\partial \hat{h}'(x,y,z,w,\lambda)$ is an non-empty convex set which could be computed by:
\begin{equation}\label{eq:dhhat'}
\begin{array}{rl}
\partial \hat{h}'(x,y,z,w,\lambda)  & = \partial (\hat{g}' - f') (x,y,z,w,\lambda)\\
&=\begin{bmatrix}
(\rho_1+\rho_2)x + 2\lambda y - 2\lambda^2 x - u \\
(\rho_1+\rho_2-2\lambda^2)y + 2\lambda(x+z) \\
\rho_1z + 2\lambda y \\
-v \\
(\rho_1+\rho_2-2(\|x\|^2+\|y\|^2))\lambda + 2y^{\top}(x+z)
\end{bmatrix}.
\end{array}
\end{equation}
with $u$ and $v$ defined in (\ref{eq:u&v}).
\paragraph{\textbf{Step 2 for $(\hat{P}_{DC})$, $(P'_{DC})$ and $(\hat{P}'_{DC})$}} As stated before, step 2 of Algorithm \ref{alg:dcaforpdc} requires solving a convex optimization problem of the form (\ref{prob:dcaforpdc}). Let us write this problem as follows: \begin{equation}
\nonumber (P^k)\qquad u^{k+1}\in \argmin\{g(u)-\langle u,v^k\rangle: u\in D\}
\end{equation}
where $u$, $u^k$, $v^k$ and $D$ are defined in (\ref{prob:dcaforpdc}).
The definitions of these vectors and set $D$ for applying Algorithm \ref{alg:dcaforpdc} to $(\hat{P}_{DC})$, $(P'_{DC})$ $(\hat{P}'_{DC})$ are given in Table \ref{tab:diffofdcas} below:

\begin{table}[ht!]
	\caption{Step 2 of DCA for $(\hat{P}_{DC})$, $(P'_{DC})$ and $(\hat{P}'_{DC})$}
	\label{tab:diffofdcas}
	\begin{center}
		\scalebox{0.78}{
			\begin{tabular}{l|c|c|c}
				$(P^k)$ & $(\hat{P}_{DC})$ & $(P'_{DC})$ & $(\hat{P}'_{DC})$ \\ 
				\hline
				$g$ & $\hat{g}$ & $g'$ & $\hat{g}'$ \\ 
				$u$ & $(x,y,z,w,\lambda)$ & $(x,y,z,w,\lambda)$ & $(x,y,z,w,\lambda)$ \\ 
				$u^{k+1}$ & $(x^{k+1},y^{k+1},z^{k+1},w^{k+1},\lambda^{k+1})$ & $(x^{k+1},y^{k+1},z^{k+1},w^{k+1},\lambda^{k+1})$ & $(x^{k+1},y^{k+1},z^{k+1},w^{k+1},\lambda^{k+1})$ \\ 
				$v^k$ & $\nabla \hat{h}(x^{k},y^{k},z^{k},w^{k},\lambda^{k})$ & $\partial h'(x^{k},y^{k},z^{k},w^{k},\lambda^{k})$ & $\partial \hat{h}'(x^{k},y^{k},z^{k},w^{k},\lambda^{k})$ \\ 
				$D$ & $\mathcal{\hat{C}}$ & $\mathcal{C}$ & $\mathcal{\hat{C}}$\\
				\hline 
			\end{tabular} 
		}
	\end{center}
\end{table}
\subsection{Solution of convex optimization problem in Step 2 of DCA} \label{subsec:socpformulation}
The problem $(P^k)$ to be solved in step 2 is a convex polynomial optimization problem for $(P_{DC})$ and $(\hat{P}_{DC})$, and a convex quadratic program for $(P'_{DC})$ and $(\hat{P}'_{DC})$. Since KKT conditions are necessary and sufficient global optimality conditions for convex optimization, thus these problems can be solved efficiently via polynomial time algorithms such as Interior-Point Methods (IPM) \cite{Boyd}. In practice, we suggest to use GUROBI \cite{Gurobi}, CPLEX \cite{Cplex} and IPOPT \cite{Ipopt} for solving the convex quadratic programs required for DCA when applied to $(P'_{DC})$ and $(\hat{P}'_{DC})$. For the convex polynomial cases in $(P_{DC})$ and $(\hat{P}_{DC})$, we propose to reformulate them as a convex quadratic program with convex quadratic constraints. In fact, our experience has shown that solving a high-order convex polynomial optimization with IPM (such as using IPOPT \cite{Ipopt} or MATLAB \verb|fmincon| \cite{Matlab}) is much more difficult in practice than solving a convex quadratic program with convex quadratic constraints. Next, we show how to make such a reformulation.

Since the function $g$ in $(P_{DC})$ has the form 
\[\begin{array}{rl}
g(x,y,z,w,\lambda) = &\|y\|^2 + \|z\|^2 + \frac{\|x+w\|^2}{4} + \frac{(\lambda^2 + \|x\|^2)^2 + (\lambda^2 + \|y\|^2)^2}{2} \\
& + \frac{(4\lambda^2+4 + \|y+x\|^2 + \|y+z\|^2)^2 + (4(\lambda+1)^2 + \|y-x\|^2 + \|y-z\|^2)^2}{32}
\end{array}\]
then by introducing an additional vector $t\in \R^8$ and the convex quadratic constraints:
\begin{equation}\label{eq:addcons}
\left\lbrace 
\begin{array}{l}
\|x\|^2\leq t_1,\\
\|y\|^2\leq t_2,\\
\|x+y\|^2\leq t_3,\\
\|y+z\|^2\leq t_4,\\
\|y-x\|^2\leq t_5,\\
\|y-z\|^2\leq t_6,\\
\lambda^2\leq t_7,\\
(\lambda+1)^2\leq t_8,\\
\end{array}
\right.
\end{equation} 
we obtain the following equivalent convex quadratic program with convex quadratic constraints
\begin{equation}
\begin{array}{l}
(x^{k+1},y^{k+1},z^{k+1},w^{k+1},\lambda^{k+1}) \in \argmin\{\bar{g}(x,y,z,w,\lambda,t) \\
-\langle (x,y,z,w,\lambda),\nabla h(x^k,y^k,z^k,w^k,\lambda^k)\rangle: (x,y,z,w,\lambda)\in \mathcal{C}, (\ref{eq:addcons}) \},
\end{array}
\end{equation}
where 
\begin{equation}\label{eq:gbar}
\begin{array}{rl}
\bar{g}(x,y,z,w,\lambda,t) = & t_2 + \|z\|^2 + \|x+w\|^2/4 + (t_7+t_1)^2/2 + (t_7+t_2)^2/2 +\\
&(4t_7 + 4 + t_3 + t_4)^2/32 + (4t_8 + t_5 + t_6)^2/32.
\end{array}
\end{equation}
This optimization problem can be solved efficiently by an IPM code such as CPLEX, GUROBI or IPOPT.
\subsection{Initial point estimation}\label{subsec:startingpoint}
Finding a good initial point for DCA is an open question which is highly depending on the specific structure of the dc program. For the case where $A\in$ PD and QEiCP satisfies the co-hyperbolic condition, we propose the following procedure to compute a potentially good initial point: 
\begin{equation}\label{eq:init1}
\left\lbrace 
\begin{array}{l}
x^0 \in \argmin \{x^{\top}Ax: e^{\top}x=1,x\geq 0\},\\
\lambda^0 = (-x^{0^{\top}}Bx^0 \pm \sqrt{(x^{0^{\top}}Bx^0)^2 - 4(x^{0^{\top}}Ax^0)(x^{0^{\top}}Cx^0)})/(2x^{0^{\top}}Ax^0),\\
y^0 = \lambda^0 x^0,\\
z^0 = \lambda^0 y^0,\\
w^0 = Az^0 + By^0 + Cx^0.
\end{array}
\right. 
\end{equation}
Since $A\in$ PD, then $x^0$ is computed by solving a strictly convex quadratic program. This problem can be efficiently solved by an IPM \cite{Boyd} or a simple extension of the Block Principal Pivoting algorithm discussed in \cite{Judice94}. 

After computing $x^0$ and due to co-hyperbolic condition, then $\lambda^0$ can be found by 
\begin{equation}\label{eq:lamba0}
\lambda^0 = \left(-x^{0^{\top}}Bx^0 \pm \sqrt{(x^{0^{\top}}Bx^0)^2 - 4(x^{0^{\top}}Ax^0)(x^{0^{\top}}Cx^0)}\right)/(2x^{0^{\top}}Ax^0).
\end{equation}
Then, by the definitions of $y,z$ and $w$, we compute $y^0 = \lambda^0 x^0, z^0 = \lambda^0 y^0$ and $w^0 = Az^0 + By^0 + Cx^0$.

Note that if $w^0\geq 0$, then this initial point is also a solution of QEiCP, and DCA starting with this point does not perform any iteration. If the co-hyperbolic condition is not satisfied, then it maybe impossible to compute $\lambda^0$ by the formula (\ref{eq:lamba0}). If in this case $C\notin S_0$, then we should choose $\lambda^0$ as an arbitrary number belonging to the interval $[l,u]$ discussed in subsection \ref{subsec:CnotinS0}.  

\subsection{Local dc decomposition}\label{subsec:localdcdecomp}
The major differences among dc programming formulations $(P_{DC})$, $(P'_{DC})$, $(\hat{P}_{DC})$ and $(\hat{P}'_{DC})$ is that $(\hat{P}_{DC})$ and $(\hat{P}'_{DC})$ are based on two parameters $\rho_1$ and $\rho_2$ given in Theorem \ref{thm:boundsofspectralradius}, while $(P_{DC})$ and $(P'_{DC})$ are based on DCSOS decompositions without a parameter. In virtue of Definition \ref{def:betterdcp}, it is easy to conclude that for problems $(\hat{P}_{DC})$ and $(\hat{P}'_{DC})$, the smaller $\rho_1$ and $\rho_2$ are, the better dc decomposition is. Moreover, since $\rho_1$ and $\rho_2$ defined in Theorem \ref{thm:boundsofspectralradius} depend on the bounds of $\lambda$, thus we can further conclude that the tighter bounds of $\lambda$ are, the better dc decomposition is. This explains why we were so interested in seeking a tighter bound for $\lambda$ in section \ref{sec:boundoflambda}.

In order to find a better dc decomposition than (\ref{eq:ghatandhhat}) and (\ref{eq:ghat'andhhat'}), we propose a \emph{local dc decomposition} strategy that is explained below.

Suppose that we have found $\lambda^k\in [l,u]$ at $k$-th iteration of DCA. Then for the $(k+1)$-th iteration of DCA, we can restricted $\lambda$ in a smaller interval such as \[\lambda \in [\lambda^k-a, \lambda^k + a] \subset [l,u]\] for some suitable parameter $a\geq 0$. For example, we can take $a=\min\{1, (\lambda^k-l)/2,(u-\lambda^k)/2 \}$. Now consider the set 
\begin{equation}
\label{eq:Chatk}
\begin{array}{rl}
\mathcal{\hat{C}}^k = &\{(x,y,z,w,\lambda)\in \mathcal{\hat{C}}: \lambda\in [\lambda^k-a, \lambda^k + a], z\in [0, p_k]^n, e^{\top}z\leq (p_k)^2, \\
& y\in [\min\{0,\lambda^k-a\},\max\{0,\lambda^k+a\}]^n \}.
\end{array}
\end{equation}
where $p_k=\max\{|\lambda^k-a|,|\lambda^k+a|\}$. We can update $\rho_1$ and $\rho_2$ on $\mathcal{\hat{C}}^k$ as follows:
\begin{equation}\label{eq:rho1kandrho2k}
\left\lbrace 
\begin{array}{l}
\rho_1^k =  2(p_k+1)^2, \\
\rho_2^k =  6p_k^2+4p_k+2.
\end{array}
\right. 
\end{equation}
Clearly, $0\leq p^k\leq p$ implies $\rho^k_1\leq \rho_1$ and $\rho^k_2\leq \rho_2$. Therefore, the dc decompositions (\ref{eq:ghatandhhat}) and (\ref{eq:ghat'andhhat'}) defined with updated parameters $\rho^k_1$ and $\rho^k_2$ on $\mathcal{\hat{C}}^k$ should be better than the ones defined with $\rho_1$ and $\rho_2$ on $\mathcal{\hat{C}}^k$. Using this idea, we propose a new DCA Algorithm \ref{alg:dcavialocaldcp} for solving $(\hat{P}_{DC})$.
\begin{algorithm}[H]
	\caption{DCA with a local dc decomposition strategy for solving $(\hat{P}_{DC})$}
	\label{alg:dcavialocaldcp}
	\begin{algorithmic}
		\STATE{\textbf{Inputs}: Given $(x^0,y^0,z^0,w^0,\lambda^0)$ with $\lambda^0\in [l,u]$, and $\epsilon_1 >0$, $\epsilon_2>0$, $\epsilon_3>0$ are the tolerances mentioned in Algorithm \ref{alg:dcaforpdc}.}
		\STATE{\textbf{Outputs}: Optimal solution $(x^*,y^*,z^*,w^*,\lambda^*)$ and optimal value $f^*$.}
		\STATE{}
		\STATE{Set $k=0$, $f^*= +\infty$, $\Delta f = +\infty$ and $\Delta X = +\infty$.}
		\WHILE{($\Delta f > \epsilon_1$ and $\Delta X > \epsilon_2$ and $f^*>\epsilon_3$)}
		\STATE{
			\begin{enumerate}
				\item[\textbf{Step 1}:] Compute $a =\min\{1, (\lambda^k-l)/2,(u-\lambda^k)/2 \}$.
				\item[\textbf{Step 2}:] Compute $\mathcal{\hat{C}}^k$ via formula (\ref{eq:Chatk}).
				\item[\textbf{Step 3}:] Compute $\hat{h}(x^k,y^k,z^k,w^k,\lambda^k)$ via formula (\ref{eq:dhhat}).
				\item[\textbf{Step 4}:] Solve convex optimization problem 
				\begin{equation}\nonumber
				\begin{array}{l}
				(x^{k+1},y^{k+1},z^{k+1},w^{k+1},\lambda^{k+1}) \in \argmin\{\hat{g}(x,y,z,w,\lambda) \\
				-\langle (x,y,z,w,\lambda),\nabla \hat{h}(x^k,y^k,z^k,w^k,\lambda^k)\rangle: (x,y,z,w,\lambda)\in \mathcal{\hat{C}}^k\}.
				\end{array}
				\end{equation}
				\item[\textbf{Step 5}:] Compute 
				$$\Delta f \leftarrow |f(x^{k+1},y^{k+1},z^{k+1},w^{k+1},\lambda^{k+1}) - f(x^{k},y^{k},z^{k},w^{k},\lambda^{k})|.$$
				$$\Delta X \leftarrow \|(x^{k+1},y^{k+1},z^{k+1},w^{k+1},\lambda^{k+1}) - (x^{k},y^{k},z^{k},w^{k},\lambda^{k})\|.$$
				$$(x^*,y^*,z^*,w^*,\lambda^*)\leftarrow (x^{k+1},y^{k+1},z^{k+1},w^{k+1},\lambda^{k+1}).$$
				$$f^* \leftarrow f(x^{*},y^{*},z^{*},w^{*},\lambda^{*}).$$
				
				Update $k\leftarrow k+1.$
		\end{enumerate}}
		\ENDWHILE
		\RETURN $(x^*,y^*,z^*,w^*,\lambda^*)$ and $f^*$.
	\end{algorithmic}
\end{algorithm}

Accordingly, for $(\hat{P}'_{DC})$, we can replace $\hat{g}$ by $\hat{g}'$, $\nabla\hat{g}$ by $\partial \hat{g}'$, and $f$ by $f'$ in Algorithm \ref{alg:dcavialocaldcp} to get a DCA with local dc decomposition.

\begin{rmk}
	Algorithm \ref{alg:dcavialocaldcp} has the same convergence properties of Algorithm \ref{alg:dcaforpdc} stated in Theorem \ref{thm:convergenceofdca}. Moreover, since Algorithm \ref{alg:dcavialocaldcp} has a better dc decomposition than Algorithm \ref{alg:dcaforpdc} for $(\hat{P}_{DC})$ and $(\hat{P}'_{DC})$. So, it is expected Algorithm \ref{alg:dcavialocaldcp} to have a better numerical performance, which will be confirmed by the numerical results to be shown in the next section. Note that it is still hard to say whether Algorithm \ref{alg:dcaforpdc} for $(P_{DC})$ and $(P'_{DC})$ is better than Algorithm \ref{alg:dcavialocaldcp} for $(\hat{P}_{DC})$ and $(\hat{P}'_{DC})$. The performances of all these algorithms will be reported in the same section.
\end{rmk}

\section{Experimental results}
\label{sec:experiments}
In this section, we report some numerical results of our algorithms implemented and tested on MATLAB 2016b \cite{Matlab}. We use Yalmip \cite{Yalmip} to build our optimization models. Yalmip provides us a convenient interface to call many solvers such as GUROBI, CPLEX and IPOPT in MATLAB for solving SOCPs.

Our tests are performed on Dell Workstation equipped with i7-6820HQ 2.70GHz CPU, 32GB RAM, 64 bits Windows 10. We have used QEiCP test problems defined in \cite{Adly11,Bras16,Fernandes14,Iusem16}. The first problem denoted by SeegerAdlyQ($3$) has been taken from \cite{Adly11}. In all the remaining problems,
the matrix $A$ is the identity matrix and the matrices $B$ and $C$ are randomly generated with elements uniformly distributed in the intervals $[0, 1]$, $[0, 10]$, and $[0, 100]$. 

\subsection{Computing bounds for $\lambda$}
Firstly, we made some experiments for computing the three bounds estimations of $\lambda$ given in Theorems \ref{thm:boundsoflambda} and \ref{thm:boundsoflambdabis} and in \cite{Fernandes14} in order to find the tightest bound in practice. The numerical results of these experiments are reported in Table \ref{tab:boundsoflambda}. 

In this table and in the remaining ones, we use the following notations:
\begin{itemize}
	\item $n$: order of the matrices $A,B$ and $C$.
	\item Bounds $1$,$2$ and $3$: bounds of $\lambda$ in Theorems \ref{thm:boundsoflambda} and \ref{thm:boundsoflambdabis} and in \cite{Fernandes14}. 
	\item $l$: lower bound for $\lambda$.
	\item $u$: upper bound for $\lambda$.
	\item IT: number of iterations of DCA.
	\item CPU: cpu time (in seconds).
\end{itemize}
\begin{table}[h!]
	\caption{Lower and upper bounds of $\lambda$ computed by three different procedures}
	\label{tab:boundsoflambda}
	\begin{center}
		\scalebox{0.85}{
			\begin{tabular}{c|c|cc|cc|cc}
				\multirow{2}{*}{PROB} & \multirow{2}{*}{$n$} & \multicolumn{2}{c|}{Bounds $1$}  & \multicolumn{2}{c|}{Bounds $2$} & \multicolumn{2}{c}{Bounds $3$} \\ 
				& & $l$ & $u$ & $l$ & $u$ & $l$ & $u$\\
				\hline
				SeegerAdlyQ(3) & $3$ & $-15.197$ & $7.697$ & $-5.313$ & $1.813$ & $-10.875$ & $5.469$ \\ 
				\hline
				Rand(0,1,05) & $5$ & $-5.705$ & $3.243$ & $-3.540$ & $2.513$ & $-4.944$ & $2.669$ \\ 
				Rand(0,1,10) & $10$ & $-10.695$ & $5.807$ & $-5.652$ & $3.746$ & $-9.345$ & $4.903$ \\ 
				Rand(0,1,20) & $20$ & $-20.953$ & $10.951$ & $-10.911$ & $6.630$ & $-19.596$ & $10.042$ \\ 
				Rand(0,1,30) & $30$ & $-30.959$ & $15.963$ & $-16.238$ & $9.549$ & $-29.585$ & $15.037$ \\ 
				Rand(0,1,40) & $40$ & $-40.971$ & $20.972$ & $-21.271$ & $12.270$ & $-39.555$ & $20.022$ \\ 
				Rand(0,1,50) & $50$ & $-50.961$ & $25.971$ & $-26.103$ & $14.918$ & $-49.273$ & $24.886$ \\ 
				\hline
				Rand(0,10,05) & $5$ & $-47.863$ & $24.205$ & $-22.840$ & $15.606$ & $-42.789$ & $21.607$ \\ 
				Rand(0,10,10) & $10$ & $-99.508$ & $50.230$ & $-56.698$ & $34.105$ & $-95.230$ & $47.858$ \\ 
				Rand(0,10,20) & $20$ & $-199.347$ & $100.162$ & $-96.851$ & $56.714$ & $-188.383$ & $94.447$ \\ 
				Rand(0,10,30) & $30$ & $-300.109$ & $150.553$ & $-150.278$ & $86.672$ & $-289.527$ & $145.014$ \\ 
				Rand(0,10,40) & $40$ & $-400.669$ & $200.829$ & $-198.833$ & $111.831$ & $-389.041$ & $194.772$ \\ 
				Rand(0,10,50) & $50$ & $-500.785$ & $250.891$ & $-253.829$ & $140.857$ & $-489.202$ & $244.850$ \\ 
				\hline
				Rand(0,100,05) & $5$ & $-485.358$ & $242.858$ & $-242.185$ & $149.519$ & $-439.463$ & $219.978$ \\ 
				Rand(0,100,10) & $10$ & $-997.611$ & $497.982$ & $-513.415$ & $305.980$ & $-930.600$ & $465.548$ \\ 
				Rand(0,100,20) & $20$ & $-1994.165$ & $997.308$ & $-962.099$ & $570.225$ & $-1863.996$ & $932.255$ \\ 
				Rand(0,100,30) & $30$ & $-2999.717$ & $1500.358$ & $-1502.535$ & $854.526$ & $-2906.336$ & $1453.417$ \\ 
				Rand(0,100,40) & $40$ & $-3999.291$ & $2000.068$ & $-2038.123$ & $1155.180$ & $-3893.380$ & $1946.940$ \\ 
				Rand(0,100,50) & $50$ & $-4999.169$ & $2500.071$ & $-2492.850$ & $1385.791$ & $-4896.833$ & $2448.667$ \\ 
				\hline 
			\end{tabular} 
		}
	\end{center}
\end{table}

We observe in Table \ref{tab:boundsoflambda} that the Bounds $2$ computed via Theorem \ref{thm:boundsoflambdabis} are the tightest bounds for $\lambda$ and are much better than Bounds $1$ and $3$. Thus, we decide to use these Bounds $2$ to perform the remaining simulations. 

\subsection{Testing DCA without local dc decomposition}
We have tested the performance (CPU time and number of iterations) of DCA for solving $(P_{DC})$, $(P'_{DC})$, $(\hat{P}_{DC})$ and $(\hat{P}'_{DC})$ without local dc decomposition. The convex quadratic subproblems were solved by GUROBI. The numerical results are summarized in Table \ref{tab:numresultsdcas} with all $\epsilon_{i\in \IntEnt{1,3}}=\epsilon=10^{-3}$ and in Table \ref{tab:numresultsdcasbis} with all $\epsilon_{i\in \IntEnt{1,3}}=\epsilon=10^{-3}$. The starting points are computed by the heuristic described in subsection \ref{subsec:startingpoint}. The values in AVG and STD stand for the average and the standard deviation of the corresponding column. 
\begin{table}[h!]
	\caption{Numerical results of DCA for $(P_{DC})$, $(P'_{DC})$, $(\hat{P}_{DC})$ and $(\hat{P}'_{DC})$ with $\epsilon=10^{-3}$}
	\label{tab:numresultsdcas}
	\begin{center}
		\scalebox{0.74}{
			\begin{tabular}{c|ccc|ccc|ccc|ccc}
				\multirow{2}{*}{PROB} & \multicolumn{3}{c|}{DCA for $(P_{DC})$}  & \multicolumn{3}{c|}{DCA for $(\hat{P}_{DC})$} & \multicolumn{3}{c|}{DCA for $(P'_{DC})$} & \multicolumn{3}{c}{DCA for $(\hat{P}'_{DC})$} \\ 
				& $\lambda$ & IT & CPU & $\lambda$ & IT & CPU & $\lambda$ & IT & CPU & $\lambda$ & IT & CPU  \\
				\hline
				SeegerAdlyQ(3) & $-0.717$ & $6$ & $0.67$ &$-0.713$ & $5$ & $0.34$ &$-0.718$ & $5$ & $0.59$ &$-0.713$ & $3$ & $0.23$ \\
				\hline
				Rand(0,1,05) & $0.703$ & $11$ & $1.34$ &$0.680$ & $8$ & $0.64$ &$0.706$ & $10$ & $1.25$ &$0.717$ & $39$ & $3.13$ \\
				Rand(0,1,10) & $0.752$ & $7$ & $1.03$ &$0.777$ & $10$ & $0.83$ &$0.753$ & $7$ & $1.05$ &$0.759$ & $37$ & $3.26$ \\
				Rand(0,1,20) & $0.937$ & $7$ & $1.57$ &$0.960$ & $2$ & $0.18$ &$0.935$ & $6$ & $1.37$ &$0.954$ & $75$ & $7.12$ \\
				Rand(0,1,30) & $0.915$ & $6$ & $2.49$ &$0.937$ & $2$ & $0.20$ &$0.927$ & $6$ & $2.42$ &$0.932$ & $137$ & $14.96$ \\
				Rand(0,1,40) & $0.954$ & $7$ & $4.41$ &$0.971$ & $2$ & $0.23$ &$0.970$ & $7$ & $4.47$ &$0.979$ & $185$ & $23.76$ \\
				Rand(0,1,50) & $0.951$ & $6$ & $6.50$ &$0.972$ & $2$ & $0.30$ &$0.976$ & $8$ & $7.96$ &$0.977$ & $209$ & $33.41$ \\
				\hline
				Rand(0,10,05) & $0.876$ & $10$ & $2.00$ &$0.918$ & $229$ & $30.93$ &$0.880$ & $10$ & $2.05$ &$0.999$ & $2$ & $0.27$ \\
				Rand(0,10,10) & $0.961$ & $17$ & $3.92$ &$0.911$ & $2$ & $0.28$ &$0.929$ & $17$ & $3.81$ &$0.910$ & $14$ & $1.96$ \\
				Rand(0,10,20) & $1.145$ & $18$ & $5.37$ &$1.180$ & $2$ & $0.29$ &$1.171$ & $16$ & $5.00$ &$1.180$ & $9$ & $1.76$ \\
				Rand(0,10,30) & $1.074$ & $24$ & $12.12$ &$1.074$ & $2$ & $0.32$ &$1.104$ & $19$ & $9.46$ &$1.074$ & $2$ & $0.34$ \\
				Rand(0,10,40) & $1.047$ & $16$ & $12.74$ &$1.087$ & $2$ & $0.42$ &$1.097$ & $12$ & $9.93$ &$1.087$ & $2$ & $0.45$ \\
				Rand(0,10,50) & $1.035$ & $18$ & $22.67$ &$1.049$ & $2$ & $0.51$ &$1.078$ & $13$ & $15.89$ &$1.049$ & $14$ & $3.65$ \\
				\hline
				Rand(0,100,05) & $1.297$ & $39$ & $14.37$ &$1.117$ & $378$ & $88.90$ &$1.203$ & $31$ & $11.56$ &$1.127$ & $2$ & $0.48$ \\
				Rand(0,100,10) & $0.848$ & $10$ & $3.91$ &$1.015$ & $33$ & $8.04$ &$0.949$ & $8$ & $3.20$ &$1.016$ & $9$ & $2.17$ \\
				Rand(0,100,20) & $1.097$ & $12$ & $5.84$ &$1.146$ & $2$ & $0.50$ &$1.105$ & $17$ & $8.58$ &$1.146$ & $2$ & $0.52$ \\
				Rand(0,100,30) & $1.043$ & $9$ & $6.32$ &$1.047$ & $2$ & $0.55$ &$1.061$ & $12$ & $8.48$ &$1.047$ & $2$ & $0.56$ \\
				Rand(0,100,40) & $0.972$ & $10$ & $9.40$ &$1.007$ & $2$ & $0.59$ &$0.990$ & $9$ & $8.86$ &$1.007$ & $2$ & $0.63$ \\
				Rand(0,100,50) & $1.071$ & $10$ & $14.03$ &$1.049$ & $2$ & $0.67$ &$1.045$ & $8$ & $11.40$ &$1.049$ & $2$ & $0.72$ \\
				\hline 
				AVG &&	$12.8$ & $6.88$ && $36.3$ & $7.09$ && $11.6$ & $6.18$ && $39.3$ & $5.23$\\ 
				STD	&&	$7.89$ & $5.74$ &&	$95.05$ & $20.48$ && $6.11$ & $4.30$ && $63.42$ & $8.83$\\
				\hline 
			\end{tabular} 
		}
	\end{center}
\end{table}

\begin{table}[h!]
	\caption{Numerical results of DCA for $(P_{DC})$, $(P'_{DC})$, $(\hat{P}_{DC})$ and $(\hat{P}'_{DC})$ with $\epsilon=10^{-4}$}
	\label{tab:numresultsdcasbis}
	\begin{center}
		\scalebox{0.69}{
			\begin{tabular}{c|ccc|ccc|ccc|ccc}
				\multirow{2}{*}{PROB} & \multicolumn{3}{c|}{DCA for $(P_{DC})$}  & \multicolumn{3}{c|}{DCA for $(\hat{P}_{DC})$} & \multicolumn{3}{c|}{DCA for $(P'_{DC})$} & \multicolumn{3}{c}{DCA for $(\hat{P}'_{DC})$} \\ 
				& $\lambda$ & IT & CPU & $\lambda$ & IT & CPU & $\lambda$ & IT & CPU & $\lambda$ & IT & CPU  \\
				\hline
				SeegerAdlyQ(3) & $-0.722$ & $13$ & $1.38$ &$-0.713$ & $5$ & $0.38$ &$-0.723$ & $11$ & $1.27$ &$-0.713$ & $3$ & $0.25$ \\
				\hline
				Rand(0,1,05) & $0.712$ & $32$ & $3.94$ &$0.720$ & $176$ & $14.66$ &$0.711$ & $32$ & $4.03$ &$0.726$ & $91$ & $7.71$ \\
				Rand(0,1,10) & $0.746$ & $21$ & $2.97$ &$0.757$ & $195$ & $17.15$ &$0.747$ & $20$ & $3.03$ &$0.759$ & $129$ & $11.84$ \\
				Rand(0,1,20) & $0.940$ & $26$ & $5.71$ &$0.951$ & $181$ & $17.15$ &$0.934$ & $18$ & $4.00$ &$0.956$ & $130$ & $13.17$ \\
				Rand(0,1,30) & $0.928$ & $40$ & $16.04$ &$0.937$ & $12$ & $1.45$ &$0.941$ & $20$ & $7.94$ &$0.935$ & $197$ & $23.04$ \\
				Rand(0,1,40) & $0.959$ & $22$ & $13.82$ &$0.971$ & $2$ & $0.24$ &$0.976$ & $16$ & $9.63$ &$0.982$ & $254$ & $33.99$ \\
				Rand(0,1,50) & $0.966$ & $27$ & $25.86$ &$0.972$ & $2$ & $0.31$ &$0.997$ & $33$ & $30.02$ &$0.984$ & $388$ & $64.56$ \\
				\hline
				Rand(0,10,05) & $0.857$ & $36$ & $7.28$ &$0.882$ & $633$ & $88.54$ &$0.854$ & $22$ & $4.56$ &$0.966$ & $1808$ & $260.74$ \\
				Rand(0,10,10) & $0.984$ & $71$ & $15.88$ &$0.894$ & $2437$ & $340.51$ &$0.957$ & $42$ & $9.15$ &$0.910$ & $14$ & $1.96$ \\
				Rand(0,10,20) & $1.186$ & $65$ & $18.75$ &$1.167$ & $1023$ & $150.53$ &$1.175$ & $40$ & $11.80$ &$1.180$ & $10$ & $1.49$ \\
				Rand(0,10,30) & $1.161$ & $112$ & $50.37$ &$1.067$ & $927$ & $151.06$ &$1.183$ & $72$ & $32.96$ &$1.074$ & $2$ & $0.36$ \\
				Rand(0,10,40) & $1.038$ & $46$ & $33.48$ &$1.087$ & $2$ & $0.46$ &$1.115$ & $32$ & $23.52$ &$1.087$ & $2$ & $0.46$ \\
				Rand(0,10,50) & $1.047$ & $63$ & $69.56$ &$1.049$ & $2$ & $0.51$ &$1.139$ & $51$ & $54.25$ &$1.049$ & $15$ & $4.02$ \\
				\hline
				Rand(0,100,05) & $1.323$ & $71$ & $24.17$ &$1.112$ & $6289$ & $1473.80$ &$1.228$ & $95$ & $32.51$ &$1.127$ & $2$ & $0.48$ \\
				Rand(0,100,10) & $0.866$ & $21$ & $7.70$ &$0.962$ & $9610$ & $2297.14$ &$0.955$ & $31$ & $11.36$ &$1.016$ & $10$ & $2.43$ \\
				Rand(0,100,20) & $1.100$ & $28$ & $12.40$ &$1.134$ & $10000$ & $2640.69$ &$1.142$ & $49$ & $23.22$ &$1.146$ & $2$ & $0.53$ \\
				Rand(0,100,30) & $1.046$ & $19$ & $12.60$ &$1.042$ & $4565$ & $1311.74$ &$1.076$ & $36$ & $22.80$ &$1.047$ & $2$ & $0.58$ \\
				Rand(0,100,40) & $0.967$ & $24$ & $21.49$ &$1.007$ & $2$ & $0.62$ &$0.992$ & $25$ & $22.33$ &$1.007$ & $2$ & $0.64$ \\
				Rand(0,100,50) & $1.073$ & $20$ & $25.52$ &$1.049$ & $2$ & $0.71$ &$1.054$ & $23$ & $28.11$ &$1.049$ & $2$ & $0.74$ \\
				\hline 
				AVG &&	$39.8$ & $19.42$ && $1898.2$ & $447.77$ && $35.2$ & $17.71$ && $161.2$ & $22.58$\\			 
				STD	&& $24.87$ & $16.60$ && $3184.12$ & $811.22$ && $20.03$ & $13.57$ && $402.06$ & $58.29$\\
				\hline 
			\end{tabular} 
		}
	\end{center}
\end{table}
Based on the set of problems that we have tested, we observe in Table \ref{tab:numresultsdcas} that all these algorithms got very similar computed eigenvalues, but their performances are quite different. The AVG in iterations and CPU time for DCSOS formulations ($(P_{DC})$ and $(P'_{DC})$) are smaller than those for universal dc formulations ($(\hat{P}_{DC})$ and $(\hat{P}'_{DC})$), and this difference becomes particularly evident in Table \ref{tab:numresultsdcasbis} when $\epsilon$ is decreased. The STD in iterations and CPU time for DCSOS formulations are less sensitive than for the universal dc decompositions with respect to both the increase of problem size and the decrease of tolerance $\epsilon$. For instance, DCA for $(\hat{P}_{DC})$ and $(\hat{P}'_{DC})$ perform very slowly with many iterations, but in some other cases such as Rand($0,10,30$) and Rand($0,100,50$), DCA for $(\hat{P}_{DC})$ and $(\hat{P}'_{DC})$ has a very fast convergence. This unstable issue for the universal dc decomposition seems much more visible for most of the test problems when the tolerance $\epsilon$ is decreased (see Table \ref{tab:numresultsdcasbis} for $\epsilon$ equal to $10^{-4}$). This sensitivity issue is probably caused by the impact of the parameters $\rho_1$ and $\rho_2$, since a bigger $\rho$ leads to a worse dc decomposition and yields a bad quality in DCA. Thus, their performances are hopefully to be improved by using local dc decompositions since $\rho_i^k\leq \rho, i=1,2$. 
\subsection{Testing DCA with local dc decompositions}
Table \ref{tab:numresultslocaldcas} reports the numerical results with local dc decompositions for $(\hat{P}_{DC})$ and $(\hat{P}'_{DC})$ with $\epsilon=10^{-4}$. We use GUROBI for solving the required convex quadratic programs with convex quadratic constraints. For the test problem SeegerAdlyQ(3), we use the set $\hat{\mathcal{C}}_2$ in the formulations $(\hat{P}_{DC})$ and $(\hat{P}'_{DC})$. The remaining test problems were solved by employing $\hat{\mathcal{C}}_1$.

\begin{table}[h!]
	\caption{Numerical results of DCA for $(\hat{P}_{DC})$ and $(\hat{P}'_{DC})$ with local dc decomposition and $\epsilon=10^{-4}$}
	\label{tab:numresultslocaldcas}
	\begin{center}
		\scalebox{0.74}{
			\begin{tabular}{c|ccc|ccc}
				\multirow{2}{*}{PROB} & \multicolumn{3}{c|}{DCA for $(\hat{P}_{DC})$}  & \multicolumn{3}{c}{DCA for $(\hat{P}'_{DC})$}  \\ 
				& $\lambda$ & IT & CPU & $\lambda$ & IT & CPU   \\
				\hline
				SeegerAdlyQ(3) & $-0.718$ & $39$ & $3.86$ &$-0.718$ & $39$ & $4.37$ \\
				\hline
				Rand(0,1,05) & $0.722$ & $74$ & $8.24$ &$0.722$ & $70$ & $8.25$ \\
				Rand(0,1,10) & $0.750$ & $57$ & $7.54$ &$0.754$ & $50$ & $7.04$ \\
				Rand(0,1,20) & $0.949$ & $48$ & $6.87$ &$0.953$ & $45$ & $6.76$ \\
				Rand(0,1,30) & $0.923$ & $46$ & $7.06$ &$0.944$ & $43$ & $7.62$ \\
				Rand(0,1,40) & $0.971$ & $62$ & $10.65$ &$0.995$ & $44$ & $9.03$ \\
				Rand(0,1,50) & $0.962$ & $46$ & $8.88$ &$1.000$ & $46$ & $10.85$ \\
				\hline
				Rand(0,10,05) & $0.834$ & $83$ & $10.71$ &$0.864$ & $85$ & $11.27$ \\
				Rand(0,10,10) & $0.985$ & $144$ & $19.82$ &$0.966$ & $143$ & $20.75$ \\
				Rand(0,10,20) & $1.141$ & $73$ & $11.12$ &$1.210$ & $109$ & $18.72$ \\
				Rand(0,10,30) & $1.164$ & $114$ & $19.83$ &$1.209$ & $135$ & $26.77$ \\
				Rand(0,10,40) & $1.066$ & $50$ & $9.66$ &$1.102$ & $62$ & $14.43$ \\
				Rand(0,10,50) & $1.069$ & $64$ & $13.48$ &$1.126$ & $116$ & $30.19$ \\
				\hline
				Rand(0,100,05) & $1.290$ & $276$ & $41.56$ &$1.251$ & $97$ & $14.76$ \\
				Rand(0,100,10) & $0.826$ & $85$ & $13.41$ &$0.954$ & $73$ & $11.96$ \\
				Rand(0,100,20) & $1.195$ & $111$ & $19.46$ &$1.225$ & $211$ & $40.44$ \\
				Rand(0,100,30) & $1.047$ & $185$ & $35.67$ &$1.115$ & $104$ & $23.05$ \\
				Rand(0,100,40) & $0.972$ & $64$ & $13.65$ &$1.042$ & $73$ & $18.25$ \\
				Rand(0,100,50) & $1.070$ & $69$ & $16.33$ &$1.070$ & $65$ & $18.23$ \\
				\hline
				AVG && $88.9$ & $14.62$ && $84.7$ & $15.93$\\ 
				STD	&& $56.96$ & $9.40$ && $43.00$ & $9.04$\\ 
				\hline 
			\end{tabular} 
		}
	\end{center}
\end{table}

By comparing these results with those reported in Table \ref{tab:numresultsdcasbis}, we observed that, in most of tested cases, the AVG for iterations and CPU time of DCA with local dc decomposition are much smaller than those for the algorithms without local dc decomposition. The STD values for local dc decomposition are also better than those obtained without local dc decomposition. Therefore, we can conclude that when local decomposition is used, the performances of DCA for universal dc decompositions $(\hat{P}_{DC})$ and $(\hat{P}'_{DC})$ can be improved. Moreover, the fastest and most stable algorithms with respect to AVG and STD values appear to be still the DCA for DCSOS formulations. This good performance for QEiCP indicates that the DCSOS decompositions for polynomial function should be considered as a promising technique to improve the performance of DCA for solving other polynomial optimization problems.  

\section{Conclusions}
\label{sec:conclusions}
In this paper, we have proposed several dc programming formulations for Quadratic Eigenvalue Complementarity Problem. A DCSOS decomposition for polynomial function, and a local dc decomposition technique should be considered as the most important contributions. The corresponding DCAs have been tested and the numerical results show a good performance for our methods. Moreover, we discuss a new insight in recognizing the most important key point for characterizing the quality of a dc decomposition, and propose a potentially good initial point for QEiCP. We believe that these contributions will be important tools for helping to answer these two most common open questions that arise in dc programming. 

\section*{Acknowledgments}
We would like to thank the editor and two anonymous reviewers for their constructive comments, which helped us to improve the manuscript. The research of Yi-Shuai Niu is funded by the National Natural Science Foundation of China (Grant No: 11601327) and by the Key Construction National ``985" Program of China (Grant No: WF220426001). The research of Joaquim J\'{u}dice was partially supported in the scope of R\&D Unit UID/EEA/5008/2013, financed by the applicable financial framework (FCT/MEC) through national funds and when applicable co-funded by FEDER-PT2002 Partnership Agreement.

\section*{References}

\bibliography{references}

\begin{thebibliography}{10}
\expandafter\ifx\csname url\endcsname\relax
  \def\url#1{\texttt{#1}}\fi
\expandafter\ifx\csname urlprefix\endcsname\relax\def\urlprefix{URL }\fi
\expandafter\ifx\csname href\endcsname\relax
  \def\href#1#2{#2} \def\path#1{#1}\fi

\bibitem{Seeger11}
A.~Seeger, Quadratic eigenvalue problems under conic constraints, SIAM Journal
  on Matrix Analysis and Applications 32 (2011) 700--721.

\bibitem{Seeger99}
A.~Seeger, Eigenvalue analysis of equilibrium processes defined by linear
  complementaritv conditions, Linear Algebra and Its Applications 294 (1999)
  1--14.

\bibitem{Judice09}
J.~J. J\'{u}dice, H.~D. Sherali, I.~Ribeiro, S.~Rosa, On the asymmetric
  eigenvalue complementarity problem, Optimization Methods and Software 24
  (2009) 549--586.

\bibitem{Bras12}
C.~P. Br\'{a}s, M.~Fukushima, J.~J. J\'{u}dice, S.~Rosa, Variational inequality
  formulation for the asymmetric eigenvalue complementarity problem and its
  solution by means of a gap function, Pacific Journal of Optimization 8 (2012)
  197--215.

\bibitem{Costa10}
A.~P.~D. Costa, A.~Seeger, Cone constrained eigenvalue problems, theory and
  algorithms, Computational Optimization and Applications 45 (2010) 25--57.

\bibitem{Costa04}
A.~P.~D. Costa, J.~A.~C. Martins, I.~N. Figueiredo, J.~J. J\'{u}dice, The
  directional instability problem in systems with frictional contacts, Computer
  Methods in Applied Mechanics and Engineering 193 (2004) 357--384.

\bibitem{Fernandes17}
R.~Fernandes, J.~J. J\'{u}dice, V.~Trevisan, Complementary eigenvalues of
  graphs, Linear Algebra \& Its Applications 527 (2017) 216--231.

\bibitem{Patrascu15}
A.~Patrascu, I.~Necoara, Efficient random coordinate descent algorithms for
  large-scale structured nonconvex optimization, Journal of Global Optimization
  61~(1) (2015) 19--46.

\bibitem{Adly15}
S.~Adly, H.~Rammal, A new method for solving second-order cone eigenvalue
  complementarity problems, Journal of Optimization Theory and Applications
  165~(2) (2015) 563--585.

\bibitem{Chen16}
Z.~M. Chen, L.~Q. Qi, A semismooth newton method for tensor eigenvalue
  complementarity problem, Computational Optimization \& Applications 65~(1)
  (2016) 109--126.

\bibitem{Costa17}
A.~P.~D. Costa, A.~Seeger, F.~M.~F. Sim$\text{\~{o}}$es, Complementarity
  eigenvalue problems for nonlinear matrix pencils, Applied Mathematics \&
  Computation 312 (2017) 134--148.

\bibitem{Fernandes16}
L.~M. Fernandes, M.~Fukushima, J.~J. J\'{u}dice, H.~D. Sherali, The
  second-order cone eigenvalue complementarity problem, Optimization Methods
  and Software 31 (2016) 24--52.

\bibitem{Zhang18}
L.~Zhang, C.~Shen, M.~Yang, J.~J. J\'{u}dice, A lanczos method for large-scale
  extreme lorentz eigenvalue problems, SIAM Journal on Matrix Analysis and
  Applications 39~(2) (2018) 611--631.

\bibitem{Yihui09}
Y.~H. Zhou, M.~S. Gowda, On the finiteness of the cone spectrum of certain
  linear transformations on euclidean jordan algebras, Linear Algebra and Its
  Applications 431~(5) (2009) 772--782.

\bibitem{Adly11}
S.~Adly, A.~Seeger, A non-smooth algorithm for cone constrained eigenvalue
  problems, Computational Optimization and Applications 49~(2) (2011) 299--318.

\bibitem{Bras17}
C.~P. Br\'{a}s, A.~Fischer, J.~J. J\'{u}dice, K.~Sch\"{o}nefeld, S.~Seifert, A
  block active set algorithm with spectral choice line search for the symmetric
  eigenvalue complementarity problem, Applied Mathematics \& Computation 294
  (2017) 36--48.

\bibitem{Facchinei}
F.~Facchinei, J.~S. Pang, Finite-dimensional Variational Inequalities and
  Complementarity Problems, Springer Science \& Business Media, 2007.

\bibitem{Fernandes14}
L.~M. Fernandes, J.~J. J\'{u}dice, H.~D. Sherali, M.~Fukushima, On the
  computation of all eigenvalues for the eigenvalue complementarity problem,
  Journal of Global Optimization 59 (2014) 307--326.

\bibitem{Iusem_inpress}
A.~N. Iusem, J.~J. J\'{u}dice, V.~Sessa, P.~Sarabando, Splitting methods for
  the eigenvalue complementarity problem, to appear in Optimization Methods and
  Software.

\bibitem{Judice08}
J.~J. J\'{u}dice, M.~Raydan, S.~Rosa, S.~Santos, On the solution of the
  symmetric complementarity problem by the spectral projected gradient method,
  Numerical Algorithms 37 (2008) 391--407.

\bibitem{Lethi12}
H.~A. LeThi, M.~Moeini, D.~T. Pham, J.~J. J\'{u}dice, A dc programming approach
  for solving the symmetric eigenvalue complementarity problem, Computational
  Optimization and Applications 51 (2012) 1097--1117.

\bibitem{Niu12}
Y.~S. Niu, H.~A. LeThi, D.~T. Pham, J.~J. J\'{u}dice, Efficient dc programming
  approaches for the asymmetric eigenvalue complementarity problem,
  Optimization Methods and Software 28 (2013) 812--829.

\bibitem{Bras16}
C.~P. Br\'{a}s, A.~N. Iusem, J.~J. J\'{u}dice, On the quadratic eigenvalue
  complementarity problem, Journal of Global Optimization 66~(2) (2016)
  153--171.

\bibitem{Cottle}
R.~Cottle, J.~S. Pang, R.~E. Stone, The Linear Complementarity Problem,
  Academic Press, New York, 1992.

\bibitem{Bras15}
C.~P. Br\'{a}s, M.~Fukushima, A.~N. Iusem, J.~J. J\'{u}dice, On the quadratic
  eigenvalue complementarity problem over a general convex cone, Applied
  Mathematics \& Computation 271 (2015) 594--608.

\bibitem{Fernandes}
L.~M. Fernandes, J.~J. J\'{u}dice, M.~Fukushima, A.~Iusem, On the symmetric
  quadratic eigenvalue complementarity problem, Optimization Methods and
  Software 29 (2014) 751--770.

\bibitem{Iusem17}
A.~N. Iusem, J.~J. J\'{u}dice, V.~Sessa, H.~Sherali, The second-order cone
  quadratic eigenvalue complementarity problem, Pacific Journal of Optimization
  13 (2017) 475--500.

\bibitem{Iusem16}
A.~N. Iusem, J.~J. J\'{u}dice, V.~Sessa, H.~D. Sherali, On the numerical
  solution of the quadratic eigenvalue complementarity problem, Numerical
  Algorithms 72 (2016) 721--747.

\bibitem{Pham98}
D.~T. Pham, H.~A. LeThi, Dc optimization algorithms for solving the trust
  region subproblem, SIAM Journal on Optimization 8 (1998) 476--507.

\bibitem{Pham18}
H.~A. LeThi, D.~T. Pham, Dc programming and dca: thirty years of developments,
  Math. Program., Ser. B (2018) 1--64.

\bibitem{Niu15}
Y.~S. Niu, J.~J. J\'{u}dice, H.~A. LeThi, D.~T. Pham, Solving the quadratic
  eigenvalue complementarity problem by dc programming, Modelling, Computation
  and Optimization in Information Systems and Management Sciences, Advances in
  Intelligent Systems and Computing 359 (2015) 203--214.

\bibitem{Niu17}
Y.~S. Niu, On difference-of-sos and difference-of-convex-sos decompositions for
  polynomials, arXiv preprint arXiv:1803.09900.

\bibitem{Ipopt}
A.~W\"{a}chter, L.~T. Biegler, On the implementation of a primal-dual interior
  point filter line search algorithm for large-scale nonlinear programming,
  Mathematical Programming 106~(1) (2006) 25--57.

\bibitem{Cplex}
IBM, \href{https://www.ibm.com/us-en/marketplace/ibm-ilog-cplex}{Ibm ilog cplex
  optimization studio 12.7}.
\newline\urlprefix\url{https://www.ibm.com/us-en/marketplace/ibm-ilog-cplex}

\bibitem{Gurobi}
G.~Optimization, \href{http://www.gurobi.com/}{Gurobi 7.0}.
\newline\urlprefix\url{http://www.gurobi.com/}

\bibitem{Pham05}
D.~T. Pham, H.~A. LeThi, The dc programming and dca revisited with dc models of
  real world nonconvex optimization problems, Annals of Operations Research 133
  (2005) 23--46.

\bibitem{Rockafellar}
R.~T. Rockafellar, Convex Analysis, Princeton University Press, Princeton,
  1970.

\bibitem{Niu10}
Y.~S. Niu, Programmation DC \& DCA en Optimisation Combinatoire et Optimisation
  Polynomiale via les Techniques de SDP, Minist\`{e}re de l'enseignement
  sup\'{e}rieur et de la recherche, Institut National Des Sciences
  Appliqu\'{e}es de Rouen, France, 2010.

\bibitem{Niu14}
Y.~S. Niu, D.~T. Pham, Dc programming approaches for bmi and qmi feasibility
  problems, Advanced Computational Methods for Knowledge Engineering:
  Proceedings of the 2nd International Conference on Computer Science, Applied
  Mathematics and Applications (ICCSAMA 2014) 282 (2014) 37--63.

\bibitem{Pham02}
D.~T. Pham, H.~A. LeThi, Dc programming. theory, algorithms, applications: The
  state of the art, First International Whorkshop on Global Constrained
  Optimization and Constraint Satisfaction, Nice.

\bibitem{Horst}
R.~Horst, P.~M. Pardalos, V.~T. Nguyen, Introduction to Global Optimization,
  2nd Edition, Springer US, 2000.

\bibitem{Ahmadi}
A.~A. Ahmadi, G.~Hall, Dc decomposition of nonconvex polynomials with algebraic
  techniques, Math. Program., Ser.B (2017) 1--26.

\bibitem{Niu11}
D.~T. Pham, Y.~S. Niu, An efficient dc programming approach for portfolio
  decision with higher moments, Computational Optimization and Applications
  50~(3) (2011) 525--554.

\bibitem{Boyd}
S.~Boyd, L.~Vandenberghe, Convex Optimization, 1st Edition, Cambridge
  University Press, 2004.

\bibitem{Matlab}
MathWorks, \href{http://www.mathworks.com/help/matlab/}{Matlab documentation}.
\newline\urlprefix\url{http://www.mathworks.com/help/matlab/}

\bibitem{Judice94}
J.~J. J\'{u}dice, F.~M. Pires, A block principal pivoting algorithm for
  large-scale strictly monotone linear complementarity problems, Computers \&
  Operations Research 21 (1994) 587--596.

\bibitem{Yalmip}
J.~L\"{o}fberg, \href{https://yalmip.github.io/}{Yalmip: A toolbox for modeling
  and optimization in matlab}, Proceedings of the CACSD Conference, Taipei,
  Taiwan (2004).
\newline\urlprefix\url{https://yalmip.github.io/}

\end{thebibliography}

\end{document}